\documentclass[11pt,leqno]{article}

\usepackage{graphicx}
\usepackage{latexsym} 
\usepackage{amsmath}
\usepackage{amssymb}
\usepackage{amsfonts}
\usepackage{enumerate}
\usepackage{theorem}

\theoremstyle{plain}
\newtheorem{teo}{Theorem}[section]
\newtheorem{lem}[teo]{Lemma}
\newtheorem{cor}[teo]{Corollary}
\newtheorem{prop}[teo]{Proposition}

{\theorembodyfont{\rmfamily}
\newtheorem{defin}[teo]{Definition}
\newtheorem{oss}[teo]{Remark}

}

\renewcommand{\eqref}[1]{\textnormal{(\ref{#1})}}

\numberwithin{equation}{section}

\newcommand{\cvd}{\hfill$\square$}
\newcommand{\proof}[1]{\noindent\textsc{Proof#1}}
\newcommand{\rmi}{\mathrm{i}}
\newcommand{\rme}{\mathrm{e}}

\title{Stable determination of sound-hard polyhedral scatterers by a minimal number of scattering measurements}

\author{Hongyu \textsc{Liu}\thanks{Department of Mathematics, Hong Kong Baptist University, Kowloon, Hong Kong, China. Email:  {\tt  hongyuliu@hkbu.edu.hk}},
Michele \textsc{Petrini}\thanks{Dipartimento di Matematica e Geoscienze, Universit\`a degli Studi di Trieste, Trieste, Italy.  Email: {\tt  mpetrini139@yahoo.it}}, Luca \textsc{Rondi}\thanks{Dipartimento di Matematica e Geoscienze, Universit\`a degli Studi di Trieste, Trieste, Italy.  Email: {\tt rondi@units.it}}, Jingni \textsc{Xiao}\thanks{ Department of Mathematics, Hong Kong Baptist University, Kowloon, Hong Kong, China. Email: {\tt xiaojn@live.com}} }

\date{}

\begin{document}

\maketitle

\setcounter{section}{0}
\setcounter{secnumdepth}{2}

\begin{abstract}

The aim of the paper is to establish optimal stability estimates for the determination
of sound-hard polyhedral scatterers in
$\mathbb{R}^N$, $N \geq 2$, by a minimal number of far-field measurements. This work is a significant and highly nontrivial extension of the stability estimates for the determination of sound-soft polyhedral scatterers by far-field measurements, proved by one of the authors, to the much more challenging sound-hard case.

The admissible polyhedral scatterers satisfy minimal a priori assumptions of
Lipschitz type and may include at the same time solid obstacles
and screen-type components. In this case we obtain a stability estimate with $N$ far-field measurements. Important features of such an estimate are that we have an explicit dependence on the parameter $h$ representing the minimal size of the cells forming the boundaries of the admissible polyhedral scatterers, and that the modulus of continuity, provided the error is small enough with respect to $h$, does not depend on $h$. If we restrict to $N=2,3$ and to polyhedral obstacles, that is to polyhedra, then we obtain stability estimates with fewer measurements, namely first with $N-1$ measurements and then with a single measurement. In this case the dependence on $h$ is not explicit anymore and the modulus of continuity depends on $h$ as well.

\medskip


\medskip

\end{abstract}

\section{Introduction}\label{intro}

A set $\Sigma \subset \mathbb{R}^N$, $N\geq 2$, is called a \emph{scatterer} if it is a compact set such that $\mathbb{R}^N \backslash \Sigma$ is connected. A scatterer is said to be an obstacle if it is the closure of an open set and it is said to be a screen if its interior is empty.

If an incident time-harmonic acoustic wave encounters a scatterer then it is perturbed through the creation of a scattered or reflected wave. The total wave is given by the superposition of the incident and the scattered wave and it is characterized by the total field $u$, solution to the following exterior boundary value problem
 $$\left\{\begin{array}{ll}
\Delta u + k^2u=0 & \text{in }\mathbb{R}^N\backslash \Sigma\\
u=u^i+u^s & \text{in }\mathbb{R}^N\backslash \Sigma\\
B.C. & \text{on }\partial \Sigma\\
\displaystyle{\lim_{r\to \infty}r^{(N-1)/2}\left(\frac{\partial u^s}{\partial r}-\rmi ku^s\right)=0} & r=\|x\|.
\end{array}\right.
$$
Here $k>0$ in the reduced wave equation, or Helmholtz equation, is the wavenumber and $u^i$ is the incident field, that is the field of the incident wave.
The incident field is usually an entire solution of the Helmholtz equation, here we shall always assume that the incident wave is a time-harmonic plane wave with direction of propagation $v\in\mathbb{S}^{N-1}$, that is $u^i(x)=\rme^{\rmi kx\cdot v}$, $x\in\mathbb{R}^N$. Instead
$u^s$ is the scattered field, that is the field of the scattered wave. The last limit is the Sommerfeld radiation condition and corresponds to the fact that
the scattered wave is radiating. Moreover it implies that
the scattered field has the following asymptotic behavior
$$
u^s(x)=\frac{\rme^{\rmi k\|x\|}}{\|x\|^{(N-1)/2}}\left\{u_{\infty}(\hat{x})
+O\left(\frac{1}{\|x\|}\right)\right\},
$$
where $\hat{x}=x/\|x\|\in \mathbb{S}^{N-1}$ and
$u_{\infty}$ is the so-called \emph{far-field pattern} of $u^s$. We shall also write $u_{\infty}(\hat{x};\Sigma,k,v)$ to specify its dependence on the observation direction $\hat{x}\in\mathbb{S}^{N-1}$, the scatterer $\Sigma$, the wavenumber $k>0$ and the direction of propagation of the incident field $v\in\mathbb{S}^{N-1}$.

Finally, the boundary condition on the boundary of $\Sigma$ depends on the physical properties of the scatterer $\Sigma$. If $\Sigma$ is \emph{sound-soft}, then $u$ satisfies a homogeneous Dirichlet condition whereas if $\Sigma$ is 
\emph{sound-hard} we have a homogeneous Neumann condition.
We remark that other conditions such as the impedance boundary condition or transmission conditions for penetrable scatterers may be of interest for the applications.

The inverse scattering problem consists of recovering the scatterer $\Sigma$ by
its corresponding far-field measurements for one or more incident waves. Such an inverse problem is of fundamental importance to many areas of science and technology including radar and sonar applications, geophysical exploration, medical imaging and nondestructive testing. For a general introduction on this inverse problem see for instance \cite{Col-Kre98,Isak06}.

Physically, a far-field measurement is obtained by sending an incident plane wave and measuring the scattered wave field faraway at every possible observation directions, namely by measuring the far-field pattern $u_{\infty}$ of $u^s$.

If we measure the far-field pattern for just one incident plane wave, then we say that we use a single far-field measurement. We can obtain multiple far-field measurements by sending different incident plane waves, changing either the wavenumber or the incident  direction of propagation, and measuring the corresponding far-field patterns. In this paper we shall assume that the wavenumber $k$ is fixed and, in order to perform more measurements, we shall modify the incident direction of propagation. 

It is readily seen that the inverse problem is nonlinear and that it is formally determined with a single far-field measurement.
Establishing the unique
determination result in this formally-determined case is a longstanding
problem in the inverse scattering theory. 

The first uniqueness
result is due to Schiffer who proved the determination of a sound-soft obstacle by infinitely many far-field measurements, see \cite{L-P}. This result was improved to the case of finitely many measurements for obstacles in \cite{Col-Sle} and for screens in \cite{Ron03}. Stability estimates for the sound-soft case was proved in \cite{Isak92,Isak93}

For what concerns the sound-hard case, following the method developed in \cite{Isak90} for the transmission conditions, uniqueness for the determination of sound-hard obstacles by infinitely many far-field measurements was shown in \cite{Kir-Kre}. The same result was also obtained in the case of the impedance boundary condition.

If one reduces to a particular class of scatterers, namely the one of polyhedral scatterers, then the number of measurements needed may be considerably reduced.
The first contribution in this direction may be found in \cite{Che-Yam} where polyhedral obstacles in dimension $2$, with a suitable further non-trapping condition, were considered.

In \cite{Ale-Ron} the uniqueness for a general sound-soft polyhedral scatterer with a single measurement was proved in any dimension $N\geq 2$.
In \cite{Liu-Zou} the uniqueness for a general sound-hard polyhedral scatterer with $N$ measurements was established, again in any dimension $N\geq 2$. It was further shown in \cite{Liu-Zou3} that the number of measurements may not be reduced if sound-hard screens are allowed. However, if one considers only polyhedral obstacles, that is polyhedra, then a single measurement is enough in any dimension $N\geq 2$. This result was proved first for $N=2$, \cite{Els-Yam1}, and then extended to any $N\geq 3$, \cite{Els-Yam2}.

Concerning stability results for the determination of polyhedral scatterers by a minimal number of far-field measurements, the only result available in the literature may be found in \cite{Ron08}, where stability estimates for the determination of sound-soft polyhedral  scatterers in $\mathbb{R}^3$ with a single measurement were established. The admissible polyhedral scatterers are there assumed to satisfy essentially minimal regularity assumptions of Lipschitz type and the stability estimate is optimal, although of a logarithmic type. A particularly interesting feature of such an estimate is
that there is an explicit dependence on the parameter $h$, $h$ representing the minimal size of the cells forming the boundaries of the admissible polyhedral scatterers, and that the modulus of continuity, provided the error is small enough with respect to $h$, does not depend on this size parameter $h$.

In this work we extend the stability results of \cite{Ron08} to the more challenging case of sound-hard polyhedral scatterers.
In order to deal with sound-hard scatterers, especially when we consider determination of polyhedra with fewer measurements, there are many highly technical modifications. Moreover, there are significant extensions with respect to the sound-soft case as considered in \cite{Ron08}
that we shall briefly discuss in what follows.

We begin by establishing the stability for the determination of sound-hard polyhedral scatterers of general type, that may include, for instance, obstacles and screens at the 
same time. We consider the general case of $\mathbb{R}^N$, with $N\geq 2$.
In this case the number of far-field measurements that are required for uniqueness, thus for stability as well, can not be reduced to a number less than $N$.
The stability result for the determination of sound-hard polyhedral scatterers in $\mathbb{R}^N$
by $N$ far-field measurements is contained in Theorem~\ref{mainteoN}.

The strategy that we utilize to establish the stability estimate of Theorem~\ref{mainteoN}  follows a similar spirit to the one used in \cite{Ron08} for sound-soft scatterers. Apart from some modifications needed to deal with the Neumann boundary condition instead of the Dirichlet one, the main significant difference is that, in the sound-hard case, the required a priori bounds on the solution of the direct scattering problem, which need to be independent on the scatterer, are much harder to prove. This key preliminary point requires to establish suitable decay estimates of the scattered fields as $\|x\|\to+\infty$ that are uniform with respect to the scatterer $\Sigma$; see Proposition~\ref{l2bound}.
This is obtained with the help of the stability result in \cite{Men-Ron} for the solution of the direct problem with respect to the variation of the scatterer $\Sigma$.

Even if the strategy is similar, still there are significant novelties and extensions here with respect to the results contained in \cite{Ron08}. One of these is the fact that we generalize the technique from $\mathbb{R}^3$ to $\mathbb{R}^N$, with $N\geq 2$. 

More importantly, we consider a much more general and versatile class of admissible polyhedral scatterers with respect to the one used in \cite{Ron08}. Such a class is characterized by essentially minimal regularity assumptions of Lipschitz type. In the preliminary Section~\ref{prelsec}, in particular in Subsection~\ref{admsubs}, we introduce and extensively discuss several classes of admissible scatterers. These classes are extremely general and may turn out to be useful on many occasions, even not linked to scattering or inverse problems, so we believe that this subsection is of independent interest.

The use of such a new improved class of polyhedral scatterers requires solving some 
technical difficulties that are illustrated in Steps I and II of the geometric construction of Section~\ref{geomconstrsec}.

A remarkable consequence of these developments is that we can also generalize the result of \cite{Ron08} to this new class of polyhedral scatterers and to any dimension $N\geq 2$; see Theorem~\ref{mainteosoft}.

Moreover, we notice the following important features of the stability estimates of Theorems~\ref{mainteoN} and \ref{mainteosoft}.
First of all, these stability estimates are optimal, the dependence on the size parameter $h$ is explicit, and the modulus of continuity, when the error is small enough with respect to $h$, does not depend on $h$.

Finally, besides far-field measurements, we can also employ near-field measurements and even the more general near-field measurements with limited aperture; see Section~\ref{prelsec}, in particular Subsection~\ref{prelsubs}.
This is actually true for all of 
our stability results, which are indeed stated
with respect to near-field measurements with a limited aperture, rather than with respect to far-field measurements. However the results of Subsection~\ref{prelsubs} easily allow to obtain the corresponding estimates with respect to far-field or near-field measurements;
see Remark~\ref{diffmeasrem}

Having established a general stability result for the determination of sound-hard polyhedral scatterers by $N$ far-field measurements, we proceed to prove stability results for the determination of polyhedral obstacles, that is polyhedra, by fewer than $N$ measurements. In this case, for technical reasons, we limit ourselves to $N=2,3$. We are able to prove a stability result with a single measurement, see Theorem~\ref{mainteo1}. The stability estimate is still of optimal type, however we lose the explicit dependence on $h$ and the modulus of continuity depends, in a rather involved way, on $h$ as well.

In order to approach the challenging technical difficulties of the proof of Theorem~\ref{mainteo1} in a slightly simplified case, we first prove a stability results for polyhedra with $N-1$ measurements, again for $N=2,3$, see
Theorem~\ref{mainteoN-1}.

We observe that the inverse sound-hard obstacle problem with a single measurement is substantially different from the sound-soft case and requires a completely new and rather difficult analysis. 
In fact, the key difficulty, as for the uniqueness issue, is to avoid, in the reflection process used in the geometric construction of Section~\ref{geomconstrsec}, the reflection in a hyperplane whose normal is orthogonal to the incident direction of propagation and with respect to which the obstacle is symmetric.
In the $N-1$ measurements case, for any obstacle actually at most one hyperplane must be avoided. Still this is not an easy task, and an ad hoc modification of the general geometric construction of Section~\ref{geomconstrsec} is required, see Subsection~\ref{N-1meassubs}. In the one measurement case in $\mathbb{R}^3$, the problem becomes even more involved. In fact there might be several planes to be avoided and further difficulties arise since we need to take into account all of them simultaneously. This is performed in Subsection~\ref{1meassubs}.

The plan of the paper is as follows. In Section~\ref{prelsec} we discuss a few preliminaries. In particular we define and study suitable classes of admissible scatterers and we present a few basic properties of the solutions to the corresponding scattering problems. In Section~\ref{resultssec} the main stability results are stated. In Section~\ref{geomconstrsec} we present the main geometric construction.
Finally, in Section~\ref{proofssec} we conclude the proofs of our stability results.

\subsubsection*{Acknowledgements} 
The work of Hongyu Liu was supported by FRG grants from Hong Kong Baptist University, Hong Kong RGC General Research Funds, 12302415 and 405513, and the NSFC grant, No. 11371115. Luca Rondi was partly supported by Universit\`a
degli Studi di Trieste through FRA 2014 and by GNAMPA, INdAM.
 
\section{Classes of admissible scatterers and preliminaries}\label{prelsec}

The integer $N\geq 2$ shall always denote the space dimension. We notice that we always omit the dependence of constants on the space dimension $N$.

For any $x\in\mathbb{R}^N$, $N\geq 2$,  we denote $x=(x',x_N)\in\mathbb{R}^{N-1}\times \mathbb{R}$ and $x=(x'',x_{N-1},x_N)\in\mathbb{R}^{N-2}\times\mathbb{R}\times\mathbb{R}$.
For any $s>0$ and any $x\in\mathbb{R}^N$, 
$B_s(x)$ denotes the ball contained in $\mathbb{R}^N$ with radius $s$ and center $x$, whereas $B'_s(x')$ denotes the ball contained in $\mathbb{R}^{N-1}$ with radius $s$ and center $x'$.
Moreover, $B_s=B_s(0)$ and $B'_s=B'_s(0)$. For any ball $B$ centered at zero we denote $B^{\pm}=B\cap\{y\in\mathbb{R}^N:\ y_N\gtrless 0\}$.
Analogously, for any hyperplane $\Pi$ in $\mathbb{R}^N$, we use the following notation. If, with respect to a suitable Cartesian coordinate system, we have $\Pi=\{y\in\mathbb{R}^N:\ y_N=0\}$ then for any $x\in \Pi$ and any $r>0$ we denote $B^{\pm}_r(x)=
B_r(x)\cap\{y_N\gtrless 0\}$. Furthermore, we denote with $T_{\Pi}$ the reflection in $\Pi$, namely in this case for any $y=(y_1,\ldots,y_{N-1},y_N)\in\mathbb{R}^N$ we have
$T_{\Pi}(y)=(y_1,\ldots,y_{N-1},-y_N)$.
Finally, for any $E\subset \mathbb{R}^N$, we denote $B_s(E)=\bigcup_{x\in E}B_s(x)$.

Given a point $x\in\mathbb{R}^N$, a vector $v\in\mathbb{S}^{N-1}$, and constants
$r>0$ and $\theta$, $0<\theta\leq\pi/2$, we call $\mathcal{C}(x,v,r,\theta)$ the open cone with vertex in $x$, bisecting vector given by $v$, radius $r$ and amplitude given by $\theta$, that is
$$\mathcal{C}(x,v,r,\theta)=\left\{y\in\mathbb{R}^N:\ 0<\|y-x\|<r\text{ and }\cos(\theta)<\frac{y-x}{\|y-x\|}\cdot v\leq 1\right\}.$$
We remark that by a cone we always mean a bounded not empty open cone of the kind defined above.

By $\mathcal{H}^s$, $0\leq s\leq N$, we denote the $s$-dimensional Hausdorff measure in $\mathbb{R}^N$. We recall that $\mathcal{H}^N$ coincides with the Lebesgue measure.

\subsection{Classes of admissible scatterers and obstacles}\label{admsubs}

We recall that by a \emph{scatterer} in $\mathbb{R}^N$ we mean a compact set $\Sigma$
contained in $\mathbb{R}^N$ such that $\mathbb{R}^N\backslash \Sigma$ is connected. We say that a scatterer $\Sigma$ is an \emph{obstacle} if $\Sigma=\overline{\Omega}$ where $\Omega$ is an open set. If the interior of $\Sigma$ is empty then we usually call it a \emph{screen}. If $\Sigma$ is a scatterer in $\mathbb{R}^N$ we shall denote $G=\mathbb{R}^N\backslash \Sigma$, which is then a connected open set containing the exterior of a ball.

A more quantitative assumption on the connectedness of $G=\mathbb{R}^N\backslash \Sigma$ is the following.
Let $\delta:(0,+\infty)\to (0,+\infty)$ be a nondecreasing left-continuous function.
Let $\Sigma$ be a compact set contained in $\mathbb{R}^N$. We say that $\Sigma$ satisfies the \emph{uniform exterior connectedness} with function $\delta$ if
for any $t>0$,
for any two points $x_1$, $x_2\in\mathbb{R}^N$ so that
$B_t(x_1)$ and $B_t(x_2)$ are contained in $\mathbb{R}^N\backslash \Sigma$, and for any $s$, $0<s<\delta(t)$,
then we can find a smooth (for instance $C^1$) curve $\gamma$ connecting
$x_1$ to $x_2$ so that $B_{s}(\gamma)$ is contained in
$\mathbb{R}^N\backslash \Sigma$ as well. 

Let us notice that such an assumption is closed under convergence in the Hausdorff distance and that $\delta(t)\leq t$ for any $t>0$.

We wish to define suitable classes of admissible scatterers. We begin with some definitions. 

Let $K$ be a compact subset of $\mathbb{R}^N$. We say that $K$ is a \emph{mildly} \emph{Lipschitz hypersurface}, with or without boundary, with positive constants $r$ and $L$ if the following holds.

For any $x\in K$ there exists a bi-Lipschitz function $\Phi_x: B_{r}(x)\to\mathbb{R}^N$ such that
\begin{enumerate}[a)]
\item\label{conditiona} for any $z_1$, $z_2\in B_{r}(x)$ we have
$$L^{-1}\|z_1-z_2\|\leq\|\Phi_x(z_1)-\Phi_x(z_2)\|\leq L\|z_1-z_2\|;$$
\item\label{conditionb} $\Phi_x(x)=0$ and
$\Phi_x(K\cap B_{r}(x))\subset \Pi=\{y\in\mathbb{R}^N:\ y_N=0\}$;
\newcounter{enumi_saved}
\setcounter{enumi_saved}{\value{enumi}}
\end{enumerate}
We say that $x\in K$ belongs to the interior of $K$ if there exists $\delta$, $0<\delta\leq r$, such that $B_{\delta}(0)\cap \Pi \subset \Phi_x(K\cap B_{r}(x))$.
Otherwise we say that $x$ belongs to the boundary of $K$. We remark that the boundary of $K$ might be empty. Further we assume that

\begin{enumerate}[a)]\setcounter{enumi}{\value{enumi_saved}}
\item\label{conditionc} for any $x$ belonging to the boundary of $K$, we have that
$$\Phi_x(K\cap B_{r}(x))=\Phi_x(B_{r}(x))\cap \Pi^+$$
where $\Pi^+=\{y\in\mathbb{R}^N:\ y_N=0,\ y_{N-1}\geq 0\}$.
\end{enumerate}

Let us notice that, by compactness, such an assumption is enough to guarantee that $\mathcal{H}^{N-1}(K)$ is bounded, hence $|K|=0$. In particular, $\mathcal{H}^{N-1}(K)$ is bounded by a constant depending on the diameter of $K$, $r$ and $L$ only. Furthermore, the boundary of $K$ has $\mathcal{H}^{N-2}$ measure bounded by a constant again depending on the diameter of $K$, $r$ and $L$ only.

Moreover, $K$ has a finite number of connected components, again bounded a constant depending on the diameter of $K$, $r$ and $L$ only, and the distance between two different connected components of $K$ is bounded from below by a positive constant depending on $r$ and $L$ only.

Let us fix a bounded open set $\Omega\subset\mathbb{R}^N$, $N\geq 2$.
We shall call $\mathcal{B}(r,L,\Omega)$ the set of $K\subset\overline{\Omega}$ such that $K$ is a mildly Lipschitz hypersurface with constants $r$ and $L$. We notice that such a set is compact with respect to the Hausdorff distance, see for instance Lemma~3.6 in \cite{Men-Ron}. We finally remark that such a class is strictly related to a similar one introduced in \cite{Gia}.

Let $K$ be a compact subset of $\mathbb{R}^N$. We say that $K$ is a \emph{Lipschitz hypersurface}, with or without boundary, with positive constants $r$ and $L$ if the following holds.

For any $x\in K$, there exists a function $\varphi:\mathbb{R}^{N-1}\to\mathbb{R}$, such that $\varphi(0)=0$ and which is Lipschitz with Lipschitz constant bounded by $L$, such that, up to a rigid change of coordinates, we have $x=0$ and
\begin{equation}\label{interiorLip}
B_r(x)\cap K\subset \{y\in B_r(x):\ y_N=\varphi(y')\}.
\end{equation}
We say that $x\in K$ belongs to the interior of $K$ if there exists $\delta$, $0<\delta\leq r$, such that $B_{\delta}(x)\cap K = \{y\in B_{\delta}(x):\ y_N=\varphi(y')\}$. Otherwise we say that $x$ belongs to the boundary of $K$. We remark that the boundary of $K$ might be empty.
For any $x$ belonging to the boundary of $K$, we assume  that there exists another function
$\varphi_1:\mathbb{R}^{N-2}\to\mathbb{R}$, such that $\varphi_1(0)=0$ and which is Lipschitz with Lipschitz constant bounded by $L$, such that, up to the previous rigid change of coordinates, we have $x=0$ and
\begin{equation}\label{boundaryLip}
B_r(x)\cap K= \{y\in B_r(x):\ y_N=\varphi(y'),\ y_{N-1}\leq\varphi_1(y'')\}.
\end{equation}
We call \eqref{interiorLip} and \eqref{boundaryLip} the $L$-\emph{Lipschitz representation} of $K$ in $B_r(x)$, where \eqref{boundaryLip} is reserved for points belonging to the boundary of $K$.

We notice that a Lipschitz hypersurface with constants $r$ and $L$ is also a mildly Lipschitz hypersurface with positive constants $\tilde{r}$ and $\tilde{L}$ depending on $r$ and $L$ only.
Furthermore, we call $\mathcal{C}=\mathcal{C}(r,L,\Omega)$
the class of Lipschitz hypersurfaces with constants $r$ and $L$ contained
in $\overline{\Omega}$. We notice that $\mathcal{C}$ is compact with respect to the Hausdorff distance, too, and that $\mathcal{C}(r,L,\Omega)\subset \mathcal{B}(\tilde{r},\tilde{L},\Omega)$.

We need the following notation. For any direction $v\in \mathbb{S}^{N-1}$, we denote by
$\hat{v}$ the couple $\hat{v}=\{v,-v\}$. We also define the following distance 
$$d(\hat{v}_1,\hat{v}_2)=\min\{\|v_1-v_2\|,\|v_1+v_2\|\} \quad\text{for any }v_1,v_2\in\mathbb{S}^{N-1}.$$

Let $K$ be a compact subset of $\mathbb{R}^N$. We say that $K$ is a \emph{strongly} \emph{Lipschitz hypersurface}, with or without boundary, with positive constants $r$ and $L$ if the following holds.

First we assume that $K$ is a Lipschitz hypersurface with constants $r$ and $L$. Then we assume the following further property.
For any $x\in K$, let $e_1(x),\ldots,e_N(x)$ be the unit vectors representing the orthonormal base of the coordinate system for which the $L$-Lipschitz representation of $K$ in $B_{r}(x)$, \eqref{interiorLip} and \eqref{boundaryLip}, holds. Then $\hat{e}_N(x)$ is a Lipschitz function of $x\in K$, with Lipschitz constant bounded by $L$, and $e_{N-1}(x)$ is a Lipschitz function of $x$, as $x$ varies in the boundary of $K$, with Lipschitz constant bounded by $L$.

The usefulness of introducing the idea of strongly Lipschitz hypersurfaces is shown in the following proposition.

%Let $\Sigma$ be a scatterer such that $K=\partial \Sigma$ is a
%strongly Lipschitz hypersurface with constants $r$ and $L$, then the conclusions of Proposition~4.2 in \cite{Men-Ron} hold, with constants $0<a\leq1\leq b$ and $h_0>0$ depending on $r$ and $L$ only. Let us notice that the assumption used in \cite[Proposition~4.2]{Men-Ron} that $K$ should be oriented is not really necessary. In the next corollary we state the following interesting consequence.

\begin{prop}\label{unifconn}
Let $\Sigma$ be a scatterer such that 
$K=\partial\Sigma$ is a strongly Lipschitz hypersurface with positive constants $r$ and $L$.

Then there exists a nondecreasing left-continuous function $\delta:(0,+\infty)\to (0,+\infty)$, depending on $r$ and $L$ only, such that $\Sigma$ satisfies the uniform exterior connectedness with function $\delta$.
\end{prop}

\proof{.} Under these assumptions, the conclusions of Proposition~4.2 in \cite{Men-Ron} hold, that is, we can find constants
$0<a\leq1\leq b$ and $h_0>0$, depending on $r$ and $L$ only, and a Lipschitz function $\tilde{d}:\mathbb{R}^N\to[0,+\infty)$ such that
$$a\, \mathrm{dist}(x,\Sigma)\leq\tilde{d}(x)\leq b\, \mathrm{dist}(x,\Sigma)\quad\text{for any }x\in\mathbb{R}^N$$
and, for any $h$, $0<h\leq h_0$, $\mathbb{R}^N\backslash \Sigma_h$ is connected, where 
$\Sigma_h=\{x\in\mathbb{R}^N:\ \tilde{d}(x)\leq h\}$.
Let us notice that the assumption used in \cite[Proposition~4.2]{Men-Ron} that $K$ should be oriented is not really necessary.

Therefore, fixed $t>0$, let $x_1$, $x_2\in\mathbb{R}^N$  be any two points so that
$B_t(x_1)$ and $B_t(x_2)$ are contained in $\mathbb{R}^N\backslash \Sigma$. Then $\tilde{d}(x_i)\geq at$ for any $i=1,2$. Provided $h=at/2\leq h_0$, then $x_i\in \mathbb{R}^N\backslash \Sigma_h$ for any $i=1,2$. Then we can find a smooth (for instance $C^1$) curve $\gamma$ connecting
$x_1$ to $x_2$ so that $\gamma$ is contained in $\mathbb{R}^N\backslash \Sigma_h$. This means that any point $x$ of $\gamma$ is such that $\mathrm{dist}(x,\Sigma)\geq \tilde{d}(x)/b>at/(2b)$.That is we can choose
\begin{equation}\label{unifconnLipschitz}
\delta(t)=\left\{
\begin{array}{ll}
at/(2b) & t\in (0,2h_0/a]\\
h_0/b &t\in [2h_0/a,+\infty)
\end{array}
\right.
\end{equation}
and the proof is concluded.\cvd

\bigskip

Our next aim is to provide sufficient conditions for a Lipschitz hypersurface to be a strongly Lipschitz hypersurface. We begin with the following lemma.

\begin{lem}\label{firststronglem}
Let us fix positive constants $r$ and $L$.
Let $K\subset \tilde{K}$ be compact subsets of $\mathbb{R}^N$ such that for any $x\in K$
there exists a function $\varphi:\mathbb{R}^{N-1}\to\mathbb{R}$, such that $\varphi(0)=0$ and which is Lipschitz with Lipschitz constant bounded by $L$, such that, up to a rigid change of coordinates, we have $x=0$ and
\begin{equation}\label{interiorLip2}
B_r(x)\cap \tilde{K}= \{y\in B_r(x):\ y_N=\varphi(y')\}.
\end{equation}

Then there exist positive constants $\tilde{r}$ and $\tilde{L}$, depending on $r$ and $L$ only, such that for any $x\in K$
there exists a function $\tilde{\varphi}:\mathbb{R}^{N-1}\to\mathbb{R}$, such that $\tilde{\varphi}(0)=0$ and which is Lipschitz with Lipschitz constant bounded by $\tilde{L}$, such that, up to a rigid change of coordinates, we have $x=0$ and
\begin{equation}\label{interiorLip3}
B_{\tilde{r}}(x)\cap \tilde{K}= \{y\in B_{\tilde{r}}(x):\ y_N=\tilde{\varphi}(y')\}
\end{equation}
and the following further property holds.
For any $x\in K$, let $e_1(x),\ldots,e_N(x)$ be the unit vectors representing the orthonormal base of the coordinate system for which the $\tilde{L}$-Lipschitz representation of $\tilde{K}$ in $B_{\tilde{r}}(x)$, \eqref{interiorLip3}, holds. Then $\hat{e}_N(x)$ is a Lipschitz function of $x\in K$, with Lipschitz constant bounded by $\tilde{L}$.
\end{lem}

\proof{.}
Let us fix $x\in K$. Locally, we can give an orientation to $\tilde{K}$ near $x$, therefore without loss of generality we can assume that, locally near $x$, $\tilde{K}$ is the boundary of a Lipschitz open set.
More precisely, we can assume there exists an open set $\Omega$ such that $\tilde{K}\cap B_r(x)\subset\partial\Omega$ and, for any $y\in \tilde{K}$ whose distance from $x$ is less than $r/2$, we have
$\tilde{K}\cap B_{r/4}(y) =\partial\Omega\cap B_{r/4}(y)$ and
$$\Omega\cap B_{r/4}(y) =\{z\in B_{r/4}(y):\ z_N<\varphi(z')\},$$ 
where $\varphi$ and the orientation depend on $y$.

Let now $y_1$ and $y_2$ be two points belonging to $B_{r/16}(x)$. Let $e^1_N$ and $e^2_N$ be the corresponding vectors for which the previous Lipschitz representation holds. Then for any $y\in B_{r/8}(x)$ we can find two open cones $\mathcal{C}_1$ and $\mathcal{C}_2$, with vertex in $y$, amplitude given by an angle $\alpha_0$, $0<\alpha_0<\pi/2$ depending on $L$ only, radius $r_0=r/16$,  and bisecting vector given by $e^1_N$ and $e^2_N$ respectively such that $\mathcal{C}_i$ does not intersect $\overline{\Omega}$ whereas the opposite cone is contained in $\Omega$ for any $i=1,2$.
First we notice that the angle between $e_N^1$ and $e_N^2$ is bounded by $\pi-2\alpha_0$.
Then we take any unit vector $\nu$ on the shorter arc of the great circle on the unit sphere passing
through $e_N^1$ and $e_N^2$.

Then there exists an absolute constant $\hat{\alpha}_0$, $0<\hat{\alpha}_0<\pi/2$, such that, provided $\alpha_0\leq \hat{\alpha}_0$, we have that
the open cone $\mathcal{C}$ with vertex
in $y$, amplitude given by the angle $\alpha_1=\alpha_0/2$, radius $r_1=(\alpha_0/3)r_0$, 
and bisecting vector $\nu$ does not intersect $\overline{\Omega}$ whereas the opposite cone is contained in $\Omega$. We call this property an interior and exterior cone condition for $\Omega$ at $y\in\partial\Omega$, with amplitude $\alpha_1$, radius $r_1$ and bisecting vector $\nu$.
The proof follows from an elementary, although lengthy, geometric construction and we omit its details.

If one performs such a construction iteratively $N$ times, one obtains
$$\alpha_N=\frac{\alpha_0}{2^N}\quad\text{and}\quad r_N=\frac{\alpha_0^N}{3^N2^{N(N-1)/2}} r_0.$$

Then we subdivide the whole $\mathbb{R}^N$ into (closed) cubes with sides of length $\tilde{r}_1$ such that their diameter is less than or equal to $r/64$.
We then consider only cubes whose intersection with $K$ is not empty. Let us fix one of these and let us call it $Q$. For any vertex $x_i$, $i=1,\ldots,2^N$, of the cube $Q$ we consider a point $\tilde{x}_i\in K\cap Q$ such that
$\mathrm{dist}(x_i,K\cap Q)=\|x_i-\tilde{x}_i\|$. Then we consider $e_N^i$ as the vector corresponding to the Lipschitz representation at the point $\tilde{x}_i$. To illustrate our construction, let us assume for simplicity that $Q=[0,\tilde{r}_1]^N$. We take the points $x_1=(0,0,\ldots, 0)$ and $x_2=(\tilde{r}_1,0,\ldots,0)$ and we construct a Lipschitz function $e_N$ on the segment connecting $x_1$ and $x_2$
such that $e_N(x_i)=e_N^i$, $i=1,2$, and that, for any $x$ in such a segment, $e_N(x)$ belongs to
the shorter arc of the great circle on the unit sphere passing
through $e_N^1$ and $e_N^2$. Clearly the Lipschitz constant of such a function $e_N$ may be bounded by a constant depending on $\tilde{r}_1$ only. Then we perform the same construction on the segment connecting $(0,\tilde{r}_1,0\ldots, 0)$ and $(1,\tilde{r}_1,0\ldots,0)$ and, then, on the segments
connecting $(t,0,0,\ldots,0)$ to $(t,\tilde{r}_1,0\ldots, 0)$, for any $t$, $0\leq t\leq \tilde{r}_1$. We iterate such a construction until we find a Lipschitz function $e_N:Q\to\mathbb{S}^{N-1}$ with Lipschitz constant bounded by a constant depending on $\tilde{r}_1$ only, with the following property. For any $y\in Q\cap K$ we have that $\Omega$ satisfies an interior and exterior cone condition at any $z\in \tilde{K}\cap B_{r/16}(y)$,
with amplitude $\alpha_N$, radius $r_N$ and bisecting vector $e_N(y)$, therefore we have a Lipschitz representation of $\tilde{K}$ at $y$ with constants $\tilde{r}$ and $\tilde{L}$ depending on $\alpha_N$ and $r_N$ only, thus on $r$ and $L$ only. Performing the same construction on any cube, the proof can be concluded.\cvd

\bigskip

Let us notice that if $K$ is oriented, then we can choose 
$e_N(x)$ itself as a Lipschitz function of $x\in K$. We also observe that if $K$ is without boundary, then it is oriented, by the Jordan-Brouwer separation theorem, and we can choose $\tilde{K}=K$. Clearly these remarks applies to any connected component of $K$.
If we limit ourselves to Lipschitz hypersurfaces without boundary then we have the following corollary.

\begin{cor}\label{LipvsstrongLip}
Let us fix positive constants $r$ and $L$.
Let $K$ be a Lipschitz hypersurface with constants $r$ and $L$ without boundary. Then there exist positive constants $\tilde{r}$ and $\tilde{L}$, depending on $r$ and $L$ only, such that $K$ is a strongly Lipschitz hypersuface with constants $\tilde{r}$ and $\tilde{L}$.
\end{cor}

We conclude this discussion on sufficient conditions for a Lipschitz hypersurface to be a  strongly Lipschitz hypersurfaces by proving the following proposition.

\begin{prop}\label{finalLiplem}
Let us fix positive constants $r$ and $L$.
Let $K$ be a Lipschitz hypersurface with constants $r$ and $L$.
For any $x\in K$, let $e_1(x),\ldots,e_N(x)$ be the unit vectors representing the orthonormal base of the coordinate system for which the $L$-Lipschitz representation of $K$ in $B_{r}(x)$, \eqref{interiorLip} and \eqref{boundaryLip}, holds.

Let us assume that $\hat{e}_N(x)$ is a Lipschitz function of $x\in K$, with Lipschitz constant bounded by $L$.

Then there there exist positive constants $\tilde{r}$ and $\tilde{L}$, depending on $r$ and $L$ only, such that $K$ is a strongly Lipschitz hypersurface with constants $\tilde{r} $ and $\tilde{L}$.
\end{prop}

\proof{.} Take two couples of orthogonal vectors $e_N^1$, $e_{N-1}^1$ and $e_N^2$, $e_{N-1}^2$ for which the $L$-Lipschitz representation holds for the same point $x$
on the boundary of $K$ in a given ball of radius $r$. We notice that $e_N^2=T(e_N^1)$, where $T$ is a rotation. Then, provided the angle between $e_N^1$ and $e_N^2$ is small enough, we have that the same Lipschitz representation holds for $e_N^2$, $e_{N-1}^2$ and $e_N^2$, $T(e_{N-1}^1)$. We then apply the arguments of Lemma~\ref{firststronglem} in $\mathbb{R}^{N-1}$ and the proof may be concluded.\cvd

\bigskip

Let us observe that a sufficient condition for the assumptions of Proposition~\ref{finalLiplem} to hold has been given in Lemma~\ref{firststronglem}.

We say that an open set $D\subset\mathbb{R}^N$ is \emph{Lipschitz} with constant $r$ and $L$ if the following assumption holds.
For any $x\in\partial D$, there exists a function $\varphi:\mathbb{R}^{N-1}\to\mathbb{R}$, such that $\varphi(0)=0$ and which is Lipschitz with Lipschitz constant bounded by $L$, such that, up to a rigid change of coordinates, we have $x=0$ and
$$B_r(x)\cap D = \{y\in B_r(x): y_N<\varphi(y')\}$$
and, consequently,
$$B_r(x)\cap \partial D = \{y\in B_r(x): y_N=\varphi(y')\}.$$

Clearly, $\partial D$ is a Lipschitz hypersurface, without boundary, with the same constants $r$ and $L$. Moreover,
we notice that $D$ and $\mathbb{R}^N\backslash \overline{D}$ satisfy a uniform cone condition,
with a cone depending on $r$ and $L$ only. We recall that, given $\mathcal{C}$ a  fixed cone in $\mathbb{R}^N$, we say that an open set $D\subset\mathbb{R}^N$ satisfies the \emph{cone condition with cone} $\mathcal{C}$ if for
every $x\in D$ there exists a cone $\mathcal{C}(x)$ with vertex in $x$ and congruent to $\mathcal{C}$ such that $\mathcal{C}(x)\subset D$.

We call $\mathcal{D}=\mathcal{D}(r,L,\Omega)$ the class of sets $\partial D$ where $D\subset \Omega$ is an open set which is Lipschitz with constants $r$ and $L$. We have that
$\mathcal{D}(r,L,\Omega)\subset \mathcal{C}(r,L,\Omega)\subset \mathcal{B}(\tilde{r},\tilde{L},\Omega)$, for some constants $\tilde{r}$, $\tilde{L}$ depending on $r$ and $L$ only. Moreover, also 
$\mathcal{D}(r,L,\Omega)$ is compact with respect to the Hausdorff distance. 

We further call $\hat{\mathcal{D}}=\hat{\mathcal{D}}(r,L,\Omega)$ the class of compact sets $\Sigma\subset\overline{\Omega}$ such that $\partial \Sigma\in \mathcal{D}(r,L,\Omega)$. Also this class is compact with respect to the Hausdorff distance.

% Again, the classes of sets $\overline{\Omega}$ and $\partial\Omega$, where $\Omega$ is an open set contained in $B_R$ which is Lipschitz with constants $r$ and $L$, are compact with respect to the Hausdorff distance.

%Let $\Sigma$ be a scatterer such that $K=\partial \Sigma$ is a Lipschitz hypersurface with constants $r$ and $L$. Then each connected component of $\Sigma$ is either the closure of a Lipschitz domain or a Lipschitz hypersurface with boundary. We recall that by domain we mean a connected open set. Obviously, the exterior of any connected component of $\Sigma$ is connected, since the exterior of $\Sigma$ is. Furthermore, 
%the numbers of connected components of $\Sigma$ and $K=\partial \Sigma$ coincide, that is if a connected component of $\Sigma$ is the closure of a Lipschitz domain then its boundary is connected. Here we have made use again  of the fact that the exterior of $\Sigma$ is connected. We may conclude that
%the exterior of $\Sigma$ is connected if and only if the exteriors of its connected components are connected.

In the following classes, introduced in \cite{Men-Ron}, we combine different (mildly) Lipschitz hypersurfaces to obtain more general and complex structures.

\begin{defin}\label{classdefin}
Let us fix positive constants $r$, $L$, and a bounded open set $\Omega$.
Let us also fix $\omega:(0,+\infty)\to (0,+\infty)$ a nondecreasing left-continuous function.

We say that a compact set $K\subset \overline{\Omega}$ belongs to the class $\tilde{\mathcal{B}}=
\tilde{\mathcal{B}}(r,L,\Omega,\omega)$, respectively $\tilde{\mathcal{C}}=
\tilde{\mathcal{C}}(r,L,\Omega,\omega)$, if it satisfies the following conditions.

\begin{enumerate}[1)]
\item\label{1} $K=\bigcup_{i=1}^M K^i$ where $K^i\in\mathcal{B}(r,L,\Omega)$, respectively $\mathcal{C}(r,L,\Omega)$, for any $i=1,\ldots,M$.
\item\label{3} For any $i\in \{1,\ldots,M\}$, and any $x\in K^i$, if its distance from the boundary of $K^i$ is $t>0$, then the distance of $x$ from the union of $K^j$, with $j\neq i$, is greater than or equal to $\omega(t)$.
\end{enumerate}

We say that a compact set $\Sigma\subset \overline{\Omega}$ belongs to the class 
$\tilde{\mathcal{B}}_1=
\tilde{\mathcal{B}}_1(r,L,\Omega,\omega)$, respectively $\tilde{\mathcal{C}}_1=
\tilde{\mathcal{C}}_1(r,L,\Omega,\omega)$, if 
$\partial\Sigma\in \tilde{\mathcal{B}}(r,L,\Omega,\omega)$, respectively $\partial\Sigma\in \tilde{\mathcal{C}}(r,L,\Omega,\omega)$.
\end{defin}

We observe that, for some constants $\tilde{r}$ and $\tilde{L}$ depending on $r$ and $L$ only,
we have $\tilde{\mathcal{C}}(r,L,\Omega,\omega)\subset \tilde{\mathcal{B}}(\tilde{r},\tilde{L},\Omega,\omega)$ and $\tilde{\mathcal{C}}_1(r,L,\Omega,\omega)\subset \tilde{\mathcal{B}}_1(\tilde{r},\tilde{L},\Omega,\omega)$.
 
Let us notice that in the previous definition the number $M$ may depend on $K$. However, there exists an integer $M_0$, depending on $r$, $L$, the diameter of $\Omega$, and $\omega$ only, such that $M\leq M_0$ for any $K\in \tilde{\mathcal{B}}$, respectively $\tilde{\mathcal{C}}$.
As before, we obtain that $\mathcal{H}^{N-1}(K)$ is bounded, hence $|K|=0$. In particular $\mathcal{H}^{N-1}(K)$ is bounded by a constant depending on $r$, $L$, the diameter of $\Omega$, and $M_0$ only. Furthermore, if we set as the boundary of $K$ the union of the boundaries of $K^i$, $i=1,\ldots,M$, then
the boundary of $K$ has $\mathcal{H}^{N-2}$ measure bounded by a constant again depending on $r$, $L$, the diameter of $\Omega$, and $M_0$ only. Finally,
the number of connected components of $\mathbb{R}^N\backslash K$ is bounded by a constant $M_1$ depending on $r$, $L$, the diameter of $\Omega$, and $\omega$ only.

Without loss of generality, we shall always assume that $\omega(t)\leq t$ for any $t>0$, and that
$\lim_{t\to+\infty}\omega(t)$ is equal to a finite real number which we call $\omega(+\infty)$.

We also remark that, by Condition~\ref{3}), we have that $K^i\cap K^j$ is contained in the intersection of the boundaries of $K^i$ and $K^j$, for any $i\neq j$. By \cite[Lemma~3.8]{Men-Ron}, we have that the classes $\tilde{\mathcal{B}}$ and $\tilde{\mathcal{C}}$ are closed, and actually compact, under convergence in the Hausdorff distance.
In the next lemma we show that this is true for the classes $\tilde{\mathcal{B}}_1$ and
$\tilde{\mathcal{C}}_1$ as well.

\begin{lem}\label{compactHausdorff0}
The classes $\tilde{\mathcal{B}}_1$ and $\tilde{\mathcal{C}}_1$ introduced in Definition~\textnormal{\ref{classdefin}}
are compact under convergence in the Hausdorff distance.

Moreover, let $\Sigma$ belong to $\tilde{\mathcal{B}}_1$, or to $\tilde{\mathcal{C}}_1$,
and $x\in\partial\Sigma$. We call $G=\mathbb{R}^N\backslash \Sigma$.
For any $r_1>0$, the number of connected components $U$ of $B_{r_1}(x)\cap G$ such that $x\in\partial U$ is bounded by a constant $M_2$ depending on $r_1$, $r$, $L$, and $\omega$ only. Finally, the number of  connected components of $B_{r_1}(x)\cap G$ intersecting $B_{r_1/2}(x)$ 
is bounded by a constant $M_3$ depending on $r_1$, $r$, $L$, and $\omega$ only.
\end{lem}

\proof{.}
We begin by proving the second part of the lemma.
It is clearly enough to consider the case in which $\Sigma\in \tilde{\mathcal{B}}_1$.
Let $U$ be a connected component of $B_{r_1}(x)\cap G$ such that $x\in\partial U$. We wish to prove that there exists $s_1>0$, depending 
on $r_1$, $r$, $L$, and $\omega$ only, and 
$y$ such that
$B_{s_1}(y)\subset
U$.

Without loss of generality we can assume that $r_1\leq \tilde{r}_1$ for some $\tilde{r}_1$ depending on $r$ and $L$ only. Let $y_0\in U$ be such that
$\|y_0-x\|\leq r_1/8$ and let $y_1\in\partial\Sigma\cap\partial U$ be such that $\|x-y_1\|\leq r_1/4$  and such that $y_1$ belongs to the interior of $K^i$ for some $i\in\{1,\ldots,M\}$,
where $\partial\Sigma=K=\bigcup_{i=1}^M K^i$ as in Condition~\ref{1}). If the distance of $y_1$ from the boundary of $K^i$ is greater than $r_1/8$, then the conclusion is immediate. Otherwise, let $y_2$ be a point in the boundary of
$K^i$ whose distance from $y_1$ is less that $r_1/8$. By the local description of $K^i$ near $y_2$, we can find a point $y_3\in K^i\cap \partial U$ such that $\|x-y_3\|\leq r_1/2$ and such that its
distance 
from the boundary of $K^i$ is at least $r_1/C$ for some constant $C\geq 8$ depending on $L$ only. Then again we can conclude.

This property immediately implies that the number of connected components $U$ of 
$B_{r_1}(x)\cap G$ such that $x\in\partial U$ is bounded by a constant $M_2$ depending on $r_1$, $r$, $L$, and $\omega$ only. Moreover, it will be crucial to prove the compactness in the Hausdorff distance.
We conclude the proof of the second part with the following argument. For any $y\in B_{r_1/2}(x)\cap G$ we call $U(y)$
the connected component of 
$B_{r_1}(x)\cap G$ containing $y$.
There exists $x(y)\in \partial U(y)\cap\partial\Sigma$, such that $\|y-x(y)\|< r_1/2$ and $\|x(y)-x\|< r_1/2$. Therefore
$V(y)$, the connected component of $B_{r_1/2}(x(y))\cap G$ containing $y$ is such that $V(y)\subset U(y)$ and $x(y)\in\partial V(y)$.

We assume that there exist points $y_n\in B_{r_1/2}(x)\cap G$, for $n=1,\ldots,n_0$, such that  
$U_n=U(y_n)$ are pairwise disjoint.
Therefore, also $V_n=V(y_n)$ are pairwise disjoint, for $n=1,\ldots,n_0$.
By the previously proved property, we have that $n_0$ is bounded by a constant $M_3$ depending on $r_1$, $r$, $L$, and $\omega$ only.

%We conclude the proof of the second part with the following argument. We argue by contradiction. For any $y\in B_{r_1}(x)\cap G$ we call $U(y)$
%the connected component of 
%$B_{r_1}(x)\cap G$ containing $y$.
%We assume that there exists a sequence $\{y_n\}_{n\in\mathbb{N}}\subset B_{r_1}(x)\cap G$ such that $\lim_ny_n=x$ and 
%$U_n=U(y_n)$, $n\in\mathbb{N}$, are pairwise disjoint.
%There exists $x_n\in \partial U_n\cap\partial\Sigma$, $n\in\mathbb{N}$, such that $\lim_nx_n=x$. Without loss of generality, we assume that
%$\|y_n-x\|\leq r_1/8$ and $\|x_n-x\|\leq r_1/8$ for any $n\in\mathbb{N}$, therefore
%$V_n$, the connected component of $B_{r_1/2}(x_n)\cap G$ containing $y_n$ is such that $V_n\subset U_n$ and $x_n\in\partial V_n$. By the previously proved property, the contradiction is evident.

About compactness in the Hausdorff distance,
it is enough to prove that the class $\tilde{\mathcal{B}}_1$ is closed.
We recall that if $\Sigma_n$ converges to $\Sigma$ in the Hausdorff distance as $n\to\infty$, and we assume that $\Sigma_n$, $n\in\mathbb{N}$, and $\Sigma$ are compact sets which are uniformly bounded, 
then
$$\Sigma=\{x\in\mathbb{R}^N:\ \text{there exists }x_n\in \Sigma_n\text{ such that }\lim_nx_n=x\}.$$

We assume that $\Sigma_n\in\tilde{\mathcal{B}}_1$ converges as $n\to\infty$ to $\Sigma$. We already know that, up to subsequences that we do not relabel, $\partial\Sigma_n\to\tilde{\Sigma}\in \tilde{\mathcal{B}}$.

It is a general fact that $\partial\Sigma\subset\tilde{\Sigma}\subset\Sigma$. Hence we just need to show that $\tilde{\Sigma}=\partial\Sigma$. By contradiction, we assume that there exists $x\in\tilde{\Sigma}\backslash \partial\Sigma$. Clearly $x$ belongs to the interior of $\Sigma$, that is for some $d>0$ we have $B_d(x)\subset\Sigma$. We can find $x_n\in\partial\Sigma_n$, $n\in\mathbb{N}$, such that $\lim_nx_n=x$. We pick $r_1=d/4$ and we assume that, for $n$ large enough,
$\|x-x_n\|<d/4$. For any $n$ large enough,
there exists $y_n$ such that
$B_{s_1}(y_n)\cap\Sigma_n=\emptyset$ and $B_{s_1}(y_n)\subset
B_{d/4}(x_n)\subset B_{d/2}(x)$. Up to a subsequence, that we do not relabel, $\lim_ny_n=y\in\overline{B_{d/2}(x)}$.
But $y$ should belong to $\Sigma$, hence there exists $\tilde{y}_n\in \Sigma_n$, $n\in\mathbb{N}$, such that $\lim_n\tilde{y}_n=y$, therefore for any $n$ large enough we have that $\|\tilde{y}_n-y_n\|<s_1$ and this is a contradiction.

The argument for the class $\tilde{\mathcal{C}}_1$ is completely analogous, and the proof is concluded.\cvd

\bigskip

Finally, we consider the following definition. We recall that $T:D\to D'$, $D$ and $D'$ being open subsets of $\mathbb{R}^N$, is said to be a bi-$W^{1,\infty}$ mapping with constant $L$ if $T$ is bijective and both $\|JT\|_{L^{\infty}(D)}$ and $\|J(T^{-1})\|_{L^{\infty}(D')}$ are bounded by $L$. Here $T^{-1}$ is the inverse of $T$ and $JT$ denotes the Jacobian matrix of $T$.

\begin{defin}\label{scatdefin}
Let us fix a bounded open set $\Omega$ and
positive constants $r$, $L$, $0<r_1< r$ and $\tilde{C}>0$.
Let us also fix $\omega:(0,+\infty)\to (0,+\infty)$, % and $\delta:(0,+\infty)\to (0,+\infty)$ two 
a nondecreasing left-continuous functions.

We call $\hat{\mathcal{B}}=\hat{\mathcal{B}}(r,L,\Omega,r_1,\tilde{C},\omega)$ and
$\hat{\mathcal{C}}=\hat{\mathcal{C}}(r,L,\Omega,r_1,\tilde{C},\omega)$
 the classes of sets satisfying the following assumptions:

\begin{enumerate}[i)]
\item\label{uniformboundedness}
any $\Sigma\in \hat{\mathcal{B}}$, respectively $\hat{\mathcal{C}}$,  is a compact set contained in $\overline{\Omega}\subset\mathbb{R}^N$ such that $\Sigma$ belongs to $\tilde{\mathcal{B}}_1(r,L,\Omega,\omega)$, respectively $\tilde{\mathcal{C}}_1(r,L,\Omega,\omega)$. We call $G=\mathbb{R}^N\backslash\Sigma$.
%\item\label{Moscocompactness}
%for any $\Sigma\in\hat{\mathcal{B}}$ we have that $\partial \Sigma\in
%\tilde{\mathcal{B}}(r,L,\Omega,\omega)$;
\item\label{pcondition} for any $x\in \partial\Sigma$ and any $U$ connected
component of $G\cap B_{r_1}(x)$, with $x\in\partial U$, we can find an open set $U'$
%Let $K=\partial\Sigma=\bigcup_{i=1}^M K^i$ as in \textnormal{\ref{1})} of Definition~\textnormal{\ref{classdefin}}.
%For any $x\in \partial\Sigma$ belonging to the boundary of $K_i$, $i\in\{1,\ldots,M\}$,
%and any $U$ connected component of $G\cap B_r(x)$ we can find an open set $U_1$,
such that 
\begin{equation}\label{changecond1}
U\subset U'\subset G,%\cap \overline{B_{r_2}(x)},
\end{equation}
and a bi-$W^{1,\infty}$ mapping $T:(-1,1)^{N-1}\times (0,1)\to U'$, with constant $\tilde{C}$, such that the following
properties hold. By the regularity of $Q=(-1,1)^{N-1}\times (0,1)$, $T$
can be actually extended up to the boundary and we have that $T:\overline{Q}\to\mathbb{R}^N$ is a Lipschitz map with Lipschitz
constant bounded by $\tilde{C}$. Furthermore, if we set
$\Gamma=[-1,1]^{N-1}\times\{0\}$, we require that 
\begin{equation}\label{changecond2}
T(0)=x\quad \text{and}\quad\partial U\cap B_{r_1}(x)\subset T(\Gamma)\subset \partial G,
\end{equation}
and that,
for any $0<s<r_1$ and any $y\in U\cap B_{r_1-s}(x)$, we have
\begin{equation}\label{changecond3}
\mathrm{dist}(T^{-1}(y), \partial Q\backslash \Gamma)\geq \omega(s).
\end{equation}
\end{enumerate}
\end{defin}

%\begin{rem}\label{propertiesrem0}
%We call $U_{x}=B_{r_{x}}(x)$, and, 
%as before, we call $\tilde{U}_{x}$ the union of the connected components of $
%U_{x}\cap G$ such that $x$ belongs to their boundaries.
%We have that the number of connected components of $\tilde{U}_{x}$
%is bounded by a constant depending on $\tilde{C}$ and $r$ only. Furthermore, there exists 
%a neighbourhood of $x$ such that the number of connected components of
% $U_{x}\cap G$ intersecting such a neighbourhood is finite. This implies the following important fact. 
%If 
%$V_{x}=(U_{x}\cap G)\backslash \tilde{U}_{x}$
%and $\hat{U}_{x}=U_{x}\backslash \overline{V_{x}}$,
%then we have that $\hat{U}_{x}$ is an open set containing $x$ and such that $\hat{U}_{x}\cap G=\tilde{U}_{x}$.
%\end{rem}

\begin{oss}\label{propertiesrem}
We notice that clearly $T(\partial Q)=\partial U'$ and  
$\partial U\cap B_{r_1}(x)\subset\partial G$.

It is pointed out that, up to suitably changing the constants $r_1$ and $\tilde{C}$ involved, Condition~\textnormal{\ref{pcondition})} is satisfied provided it holds only for points $x$ belonging to the boundaries of $K_i$, $i=1,\ldots,M$, where $\partial\Sigma=\bigcup_{i=1}^MK_i$ by Condition~\textnormal{\ref{uniformboundedness}}).

Again, we have
$\hat{\mathcal{C}}(r,L,\Omega,r_1,\tilde{C},\omega)\subset \hat{\mathcal{B}}(\tilde{r},\tilde{L},\Omega,r_1,\tilde{C},\omega)$, for some constants $\tilde{r}$ and $\tilde{L}$ depending on $r$ and $L$ only,
We also remark that,  for some constants and functions depending on $r$, $L$, and
the diameter of $\Omega$ only, we have
 $\hat{\mathcal{D}}(r,L,\Omega)\subset \hat{\mathcal{C}}(\hat{r},\hat{L},\Omega,r_1,\tilde{C},\omega)$.
\end{oss}

It is emphasised that Condition~\ref{pcondition}) is an extremely weak regularity condition and that it is satisfied by rather complex structures, see for instance the discussion on sets in $\mathbb{R}^3$ satisfying this assumption in Section 4 of \cite{Ron07}, where several examples are shown.

%The Condition~\ref{pcondition}) in Definition~\ref{scatdefin}
%is a slight modification of a corresponding one introduced in \cite{LPRX}, namely the third condition of Definition~2.6. We notice that these two conditions are actually
%equivalent, up to suitably changing the constants $r_1$ and $\tilde{C}$, thus, for the purpose of the estimates needed in Section~\ref{polystabsec}, this new version makes no difference.
%The reason for this modification is that this new version is closed under convergence in the Hausdorff distance, as it is stated in the next lemma.

The following compactness result holds true.

\begin{lem}\label{compactHausdorff}
The classes $\hat{\mathcal{B}}$ and $\hat{\mathcal{C}}$ introduced in Definition~\textnormal{\ref{scatdefin}}
are compact under convergence in the Hausdorff distance.
\end{lem}

\proof{.} The argument is the same for both classes $\hat{\mathcal{B}}$ and $\hat{\mathcal{C}}$, so we limit ourselves to the first one.
It is enough to prove that the class is closed. By Lemma~\ref{compactHausdorff0},
we just need to prove that also Condition~\ref{pcondition}) is preserved in the limit.
Let $\Sigma_n\in \hat{\mathcal{B}}$, $n\in\mathbb{N}$, be such that
$\Sigma_n\to\Sigma\in \tilde{\mathcal{B}}_1$ in the Hausdorff distance as $n\to\infty$.

Let $x\in\partial\Sigma$ and let $U$ be a connected component of $G\cap B_{r_1}(x)$ with $x\in\partial U$. Let $y$ and $s>0$ be such that
$B_s(y)\subset U$ and $\|y-x\|<r_1/2$. We also consider a continuous curve $\gamma:[0,1]\to\mathbb{R}^N$ such that $\gamma(0)=y$, $\gamma(1)=x$, and
$\gamma(t)\in U$ for any $t\in[0,1)$.
Let $\tilde{U}_n$, $n\in\mathbb{N}$, be the connected component of $G_n$ containing $y$, at least for $n$ large enough.

Let $\{t_m\}_{m\in\mathbb{N}}\subset [0,1)$ be an increasing sequence 
 such that $\|\gamma(t_m)-x\|<1/m$.
Then there exists an increasing sequence $\{n_m\}_{m\in\mathbb{N}}$ of integers such that for any $n\geq n_m$ we have $\gamma([0,t_m])\subset\tilde{U}_n$.
Since there exists $\tilde{x}_n\in \Sigma_n$ converging to $x$ as $n\to\infty$, we can conclude that
there exists $x_{n_m}\in \partial\tilde{U}_{n_m}\cap\partial\Sigma_{n_m}$ such that
$\lim_mx_{n_m}=x$. It is also not difficult to show that, for any $m$ large enough, we can find
$U_{n_m}$, a connected component of $B_{r_1}(x_{n_m})\cap G_{n_m}$, such that
$x_{n_m}\in \partial U_{n_m}$, $y\in U_{n_m}$, $\|y-x_{n_m}\|<r_1/2$,
and $U_{n_m}\subset \tilde{U}_{n_m}$.

We call $T_m:Q\to U'_{n_m}$ with $U_{n_m}\subset U'_{n_m}\subset G_{n_m}$ as in Condition~\ref{pcondition}). Clearly, up to a subsequence that we do not relabel, 
$T_m$ converges uniformly on $\overline{Q}$ to $T:\overline{Q}\to \mathbb{R}^N$,
$T$ being a Lipschitz function with constant $\tilde{C}$.
Obviously $T(0)=x$ and a
straightforward computation shows that $T|_{Q}$ is actually bi-$W^{1,\infty}$, with 
constant $\tilde{C}$, between $Q$ and $U'$.
We have that $U'$ is connected and 
we need to show that $U'\cap \Sigma=\emptyset$, that is $U'\subset G$.

We assume, by contradiction, that there exists $w\in Q$ such that $T(w)\in \Sigma$. By the bi-$W^{1,\infty}$ property, we have that $B_s(T_m(w))\subset T_m(Q)$,
for some $s>0$ independent of $m$.
There exists $y_n\in\Sigma_n$,
$n\in\mathbb{N}$, such 
that $\lim_ny_n= T(w)$. On the other hand,
$\lim_m T_m(w)=T(w)$ as well, hence for $m$ large enough we have
$\|y_{n_m}-T_m(w)\|< s$ and this is a contradiction to the fact that $T_m(Q)\cap\Sigma_{n_m}=\emptyset$.

Next, we prove the first inclusion of \eqref{changecond1}. Let $x_1\in
 U$ be fixed. There exists a continuous curve $\gamma_1$ in $U$ connecting
$x_1$ with $y$. We have that, for some $d>0$, $B_d(\gamma_1)\subset U$, therefore, for any $m$ large enough, 
$\Sigma_{n_m}\cap B_{d/2}(\gamma_1)=\emptyset$
and $B_{d/2}(\gamma_1)\subset B_{r_1}(x_{n_m})$. Therefore, 
$B_{d/2}(\gamma_1)\subset U_{n_m}$ and in particular $B_{d/2}(x_1)\subset U_{n_m}\subset U'_{n_m}$. By a reasoning completely analogous to the one used to prove that $U'\subset G$, we conclude that
$x_1\in U'$.

For what concerns \eqref{changecond2} and \eqref{changecond3}, these can be proved with straightforward modifications of the above arguments and
the proof is complete.\cvd

\bigskip

Now we are ready to define the following classes of admissible scatterers.

\begin{defin}\label{scatdefin1}
Let us fix positive constants $r$, $L$ and $R$, $0<r_1< r$ and $\tilde{C}>0$.
Let us also fix $\omega:(0,+\infty)\to (0,+\infty)$ and $\delta:(0,+\infty)\to (0,+\infty)$ two  nondecreasing left-continuous functions.

We call $\hat{\mathcal{B}}_{scat}=\hat{\mathcal{B}}_{scat}(r,L,R,r_1,\tilde{C},\omega,\delta)$  the class of compact sets $\Sigma$ such that $\Sigma$ belongs to $\hat{\mathcal{B}}(r,L,B_R,r_1,\tilde{C},\omega)$ 
and satisfies
the uniform exterior connectedness with function $\delta$.

We also define $\tilde{\mathcal{B}}_{scat}=\tilde{\mathcal{B}}_{scat}(r,L,R,\omega,\delta)$ the class of compact sets $\Sigma$ belonging to $\tilde{\mathcal{B}}_1(r,L,B_R,\omega)$ and 
satisfying
the uniform exterior connectedness with function $\delta$.

Completely analogous definitions may be given for $\hat{\mathcal{C}}_{scat}$ and 
$\tilde{\mathcal{C}}_{scat}$.

We further call $\hat{\mathcal{D}}_{obst}=\hat{\mathcal{D}}_{obst}(r,L,R)$ the class of compact sets $\Sigma$ belonging to $\hat{\mathcal{D}}(r,L,B_R)$
and such that $G=\mathbb{R}^3\backslash \Sigma$ is connected.
\end{defin}

Obviously, we have $\hat{\mathcal{B}}_{scat}(r,L,R,r_1,\tilde{C},\omega,\delta)\subset
 \tilde{\mathcal{B}}_{scat}(r,L,R,\omega,\delta)$ and the same relation holds between 
 $\hat{\mathcal{C}}_{scat}$ and $\tilde{\mathcal{C}}_{scat}$. Moreover, the same relations as before hold between the classes $\hat{\mathcal{B}}_{scat}$ and
 $\tilde{\mathcal{B}}_{scat}$ and the corresponding classes $\hat{\mathcal{C}}_{scat}$ and $\tilde{\mathcal{C}}_{scat}$.
We notice that any scatterer $\Sigma\in \hat{\mathcal{D}}_{obst}$ is indeed an obstacle, that is,
$\Sigma$ is the closure of its interior which is a bounded open set with Lipschitz boundary, with constants $r$ and $L$.
By Corollary~\ref{LipvsstrongLip} and Proposition~\ref{unifconn}, 
for some constants and functions depending on $r$, $L$, and $R$ only, we have
 $\hat{\mathcal{D}}_{obst}(r,L,R)\subset \hat{\mathcal{C}}_{scat}(\tilde{r},\tilde{L},R,r_1,\tilde{C},\omega,\delta)$. 

By our earlier discussion, in particular by Lemmas~\ref{compactHausdorff0} and \ref{compactHausdorff}, it is easy to note that all these classes
$\tilde{\mathcal{B}}_{scat}$, $\hat{\mathcal{B}}_{scat}$, $\tilde{\mathcal{C}}_{scat}$, $\hat{\mathcal{C}}_{scat}$, 
and $\hat{\mathcal{D}}_{obst}$ are compact with respect to the Hausdorff distance.

Finally, the sets belonging to the class $\hat{\mathcal{B}}$, thus in particular 
scatterers belonging to the class $\hat{\mathcal{B}}_{scat}$, satisfy the following property.

\begin{prop}\label{suffcond}
Let us fix positive constants $r$, $L$ and $R$, $0<r_1<r$, and $\tilde{C}>0$.
Let us also fix $\omega:(0,+\infty)\to (0,+\infty)$ a nondecreasing left-continuous functions.

 %and $\delta:(0,+\infty)\to (0,+\infty)$ two nondecreasing left-continuous functions.

Let $\hat{\mathcal{B}}=\hat{\mathcal{B}}(r,L,B_R,r_1,\tilde{C},\omega)$.
Then there exist constants $p>2$ and $\tilde{C}_1>0$, depending on $\hat{\mathcal{B}}$ only, such that, for any $\Sigma\in\hat{\mathcal{B}}$,
we have
\begin{equation}\label{Sobolevcondition}
\|v\|_{L^p(B_{R+1}\backslash \Sigma)}\leq \tilde{C}_1\|v\|_{H^1(B_{R+1}\backslash \Sigma)}\quad\text{for any }v\in H^1(B_{R+1}\backslash \Sigma).
\end{equation}

Moreover, the immersion of $H^1(B_{R+1}\backslash \Sigma)$ into $L^2(B_{R+1}\backslash\Sigma)$ is compact, for any $\Sigma\in\hat{\mathcal{B}}$.
\end{prop}

%We immediately notice that also the property described in \eqref{Sobolevcondition} is closed with respect to the Hausdorff convergence. ??? vero ???

\proof{.} We fix $\Sigma\in\hat{\mathcal{B}}$ and we call $G=\mathbb{R}^3\backslash \Sigma$. Let us take $v$ belonging to $H^1(B_{R+1}\backslash \Sigma)$.
Without loss of generality, by an easy extension argument around $\partial B_{R+1}$, we can assume that $v$ actually belongs to $H^1(\mathbb{R}^3\backslash \Sigma)$, it has bounded support, and its $H^1$ norm is controlled by a constant $C$, depending on $R$ only, times
 the corresponding $H^1$ norm in $B_{R+1}\backslash \Sigma$.

We start with a local construction.
We fix $x\in \partial \Sigma$ and $U$ a connected component of
$B_{r_1}(x)\cap G$ such that $x\in\partial U$. We consider $U'$
and $T:Q\to U'$ as in Condition~\ref{pcondition}) of Definition~\ref{scatdefin}.

Clearly we have that $U'$ satisfies
\begin{equation}\label{Sobolevprel2}
\|v\|_{L^{s_1}(U')}\leq C_1\|v\|_{H^1(U')}\quad\text{for any }v\in H^1(U')
\end{equation}
for some constants $s_1>2$ and $C_1>0$.
Since $U\cap B_{3r_1/4}(x)\subset U'$, we conclude that
\begin{equation}\label{cutoffestimatefinalbis}
\|v\|_{L^{s_1}(U\cap B_{3r_1/4}(x))}
\leq
C_1
\|v\|_{H^1(U')}\quad\text{for any }v\in H^1(U').
\end{equation}

We now consider a covering argument as follows.
For any $x\in\partial\Sigma$, let $W_n$, $n=1,\ldots,n_0$, be the connected components of $B_{r_1/2}(x)\cap G$ such that $W_n\cap B_{r_1/4}(x)\neq\emptyset$. By Lemma~\ref{compactHausdorff0}, $n_0\leq M_3$, where $M_3$ is a constant depending on $r_1$, $r$, $L$, and $\omega$ only.
Let $y_n\in W_n\cap B_{r_1/4}(x)$,
$n=1,\ldots,n_0$. As in the proof of Lemma~\ref{compactHausdorff0},
there exists $x_n\in \partial W_n\cap\partial\Sigma$, such that $\|y_n-x_n\|< r_1/4$ and $\|x_n-x\|< r_1/4$. 
We call $U_n$ the connected component of $B_{r_1}(x_n)\cap G$ containing $y_n$ and we observe that $x_n\in\partial U_n$ %, $V_n\subset W_n$, 
and $W_n\subset U_n$. Actually, $W_n\subset U_n\cap B_{3r_1/4}(x_n)$.

We conclude that for any $x\in\partial\Sigma$, there exist $n_0$ points $x_1,\ldots,x_{n_0}$, with $n_0\leq M_3$, with the following property.
For any $n=1,\ldots,n_0$, there exists
$U_n$, a connected component of
$B_{r_1}(x_n)\cap G$, such that $x_n\in\partial U_n$, and
moreover
$$B_{r_1/4}(x)\cap G\subset \bigcup_{n=1}^{n_0}(U_n\cap B_{3r_1/4}(x_n)).$$

We fix $\delta=r_1/16$ and define the compact set $A_1=\overline{B_{\delta}(\partial\Sigma)\cap G}$. We notice that
$$A_1\subset\bigcup_{x\in \partial\Sigma}B_{r_1/4}(x).$$
By the compactness of $A_1$,
we can find a finite number of points $z_i\in \partial\Sigma$, $i=1,\ldots,m_1$, such that
$$A_1\subset \bigcup_{i=1}^{m_1}B_{r_1/4}(z_i).$$
With a rather simple construction, it is possible to choose $m_1$ depending on $r_1$ and $R$ only, for instance by taking points such that $B_{r_1/16}(z_i)\cap B_{r_1/16}(z_j)$ is empty for $i\neq j$.

We further find a finite number of points $z_i\in \partial B_{R+1}$, $i=m_1+1,\ldots,m_1+m_2$, such that
$$A_2=\overline{B_{1/16}(\partial B_{R+1})}\subset\bigcup_{i=m_1+1}^{m_1+m_2}B_{1/4}(z_i),$$ with $m_2$ depending on $R$ only.

Finally, we call $r_3=\min\{1,r_1\}$ and
$$A_3=\{x\in B_{R+1}\backslash \Sigma:\ \mathrm{dist}(x,\partial(B_{R+1}\backslash \Sigma)) \geq r_3/16\}.$$
We can find points $z_i\in A_3$, $i=m_1+m_2+1,\ldots,m_1+m_2+m_3$, such that
$$A_3\subset \bigcup_{i=m_1+m_2+1}^{m_1+m_2+m_3}B_{r_3/32}(z_i).$$
Again the number $m_3$ may be bounded by a constant depending on $r_1$ and $R$ only.

By applying \eqref{cutoffestimatefinalbis}, at most $M_3$ times for any $z_i$, $i=1,\ldots,m_1$, we have
\begin{equation}
\|v\|_{L^{s_1}(A_1\cap G)}\leq C_1(M_3m_1)C\|v\|_{H^1(B_{R+1}\backslash \Sigma)}.
\end{equation}

By a completely analogous argument, we can find $s_2>2$ and $C_2$
such that
\begin{equation}
\|v\|_{L^{s_2}(A_2\cap B_{R+1})}\leq C_2m_2\|v\|_{H^1(B_{R+1}\backslash \Sigma)}.
\end{equation}

Applying a classical Sobolev inequality to $D=B_{r_3/32}(z_i)$, for $i=m_1+m_2+1,\ldots,m_1+m_2+m_3$,
we can finally find
$s_3>2$ and $C_3$
such that
\begin{equation}
\|v\|_{L^{s_3}(A_3)}\leq C_3m_3\|v\|_{H^1(B_{R+1}\backslash \Sigma)}.
\end{equation}

Picking $p=\min\{s_1,s_2,s_3\}$ we obtain that
\begin{equation}
\|v\|_{L^{p}(B_{R+1}\backslash \Sigma)}\leq \tilde{C}_1
\|v\|_{H^1(B_{R+1}\backslash \Sigma)}.
\end{equation}
It is an easy remark that $p$ and $\tilde{C}_1$ have the dependence required.

The fact that the immersion of $H^1(B_{R+1}\backslash \Sigma)$ into $L^2(B_{R+1}\backslash\Sigma)$ is compact is an immediate consequence of the property described in \eqref{Sobolevcondition}.\cvd

\bigskip

We conclude this subsection by introducing suitable classes of polyhedral scatterers.
We define a \emph{cell} as the closure of an open subset of an
$(N-1)$-dimensional hyperplane.
We say that a scatterer $\Sigma$ is \emph{polyhedral} if the boundary of $\Sigma$ is given by a finite union of cells
$\mathcal{C}_j$, $j=1,\ldots,M_1$.

Fixed positive constants $h$ and $L$,
we say that a scatterer $\Sigma$ is \emph{polyhedral with constants} $h$ \emph{and} $L$ if
the boundary of $\Sigma$ is given by a finite union of cells
$\mathcal{C}_j$, $j=1,\ldots,M_1$, where
 each $\mathcal{C}_j$ is the closure of a Lipschitz domain with constants $h$ and $L$ contained in an $(N-1)$-dimensional hyperplane and the cells are pairwise internally disjoint, that is two different cells may intersect only at boundary points.

Let $\hat{\mathcal{B}}_{scat}=\hat{\mathcal{B}}_{scat}(r,L,R,r_1,\tilde{C},\omega,\delta)$ be the class of scatterers defined in Definition~\ref{scatdefin1}.
Fixed the size parameter $h>0$, let 
$\hat{\mathcal{B}}_{scat}^h=\hat{\mathcal{B}}_{scat}^h(r,L,R,r_1,\tilde{C},\omega,\delta)$
be the set of scatterers $\Sigma\in \hat{\mathcal{B}}_{scat}$ such that $\Sigma$ is polyhedral with constants $h$ and $L$.

Analogously, let
$\hat{\mathcal{D}}_{obst}=\hat{\mathcal{D}}_{obst}(r,L,R)$ be the class of obstacles defined
in Definition~\ref{scatdefin1}. Fixed the size parameter $h>0$, let $\hat{\mathcal{D}}_{obst}^h=\hat{\mathcal{D}}_{obst}^h(r,L,R)$ be the set of obstacles $\Sigma\in \hat{\mathcal{D}}_{obst}$ such that $\Sigma$ is polyhedral with constants $h$ and $L$. Notice that in this case any $\Sigma\in \hat{\mathcal{D}}_{obst}^h$ is formed by a finite number of polyhedra.

\subsection{Preliminaries}\label{prelsubs}

In this subsection we fix positive constants $r$, $L$ and $R$, $0<r_1< r$ and $\tilde{C}>0$, and two  nondecreasing left-continuous functions
$\omega:(0,+\infty)\to (0,+\infty)$ and $\delta:(0,+\infty)\to (0,+\infty)$.
The class of admissible scatterers that we consider will be called $\mathcal{A}$. Here we pick $\mathcal{A}=\hat{\mathcal{B}}_{scat}(r,L,R,r_1,\tilde{C},\omega,\delta)$, as in Definition~\ref{scatdefin1}, and we take $\Sigma$ and $\Sigma'$ belonging to $\mathcal{A}$.

We set
\begin{equation}\label{ddefin}
d=\max\left\{\sup_{x\in\partial \Sigma\backslash\Sigma'}\mathrm{dist}(x,\partial\Sigma'),
\sup_{x\in\partial \Sigma'\backslash\Sigma}\mathrm{dist}(x,\partial\Sigma)\right\}
\end{equation}
and
\begin{equation}\label{dddefin}
\hat{d}=d_H(\partial\Sigma,\partial\Sigma')\quad\text{and}\quad\tilde{d}=d_H(\Sigma,\Sigma').
\end{equation}
We recall that $d_H$ denotes the Hausdorff distance. We notice that the following relationships among $d$, $\hat{d}$ and $\tilde{d}$ holds. First, $d$, $\hat{d}$ and $\tilde{d}$ are all bounded by $2R$. We also obviously have
$d\leq \hat{d}$.
Up to swapping the role of $\Sigma$ and $\Sigma'$, let $x\in\Sigma'$ be such that $\mathrm{dist}(x,\Sigma)=\tilde{d}$.
Clearly, $\mathrm{dist}(x,\partial\Sigma)=\tilde{d}$ as well. If $x\in \partial\Sigma'$, then we immediately conclude that $\tilde{d}\leq d$. If $x$ does not belong to $\partial\Sigma'$, then, by using the uniform exterior connectedness property of $\Sigma$, for any $s<\delta(\tilde{d})$
we can find a point $x_1\in\partial\Sigma'$ such that $\mathrm{dist}(x_1,\Sigma)=\mathrm{dist}(x_1,\partial\Sigma)\geq s$. Therefore
\begin{equation}\label{dist0}
\delta(\tilde{d})\leq d\leq \hat{d}
\end{equation}
or, in other words,
\begin{equation}\label{distances}
\tilde{d}\leq \delta^{-1}(d)\leq \delta^{-1}(\hat{d})
\end{equation}
where $\delta^{-1}:(0,+\infty)\to(0,+\infty)$ is a nondecreasing right-continuous function defined as follows
\begin{equation}\label{delta-1def}
\delta^{-1}(t)=\min\{\sup\{s:\ \delta(s)\leq t\},2R\} \quad\text{for any }t>0.
\end{equation}

On the other hand, let $x\in \partial\Sigma'$ be such that $\mathrm{dist}(x,\partial\Sigma)=\hat{d}$. If $x$ does not belong to $\Sigma$, then $\mathrm{dist}(x,\Sigma)=\hat{d}$ hence $d=\hat{d}\leq \tilde{d}$. If $x\in\Sigma$, then $B_{\hat{d}}(x)\subset\Sigma$. Hence, by the properties of the boundary of $\Sigma'$, there exists a positive constant $C_1$, depending on the class $\mathcal{A}$ only, and a point $x_1$ such that $B_{C_1\hat{d}}(x_1)\subset B_{\hat{d}}(x)\backslash\Sigma'$. We can conclude that
\begin{equation}\label{distancesopposite}
C_1 d\leq C_1\hat{d}\leq \tilde{d}\leq \delta^{-1}(d)\leq \delta^{-1}(\hat{d}).
\end{equation}

Let us notice that we also have the following property that
will be of use later on. If $C_2=(C_1+1)/C_1$, then
\begin{equation}\label{boundary}
\Sigma\triangle\Sigma'\subset \overline{B_{C_2\tilde{d}}(\partial\Sigma)}\cap \overline{B_{C_2\tilde{d}}(\partial\Sigma')}.
\end{equation}
In fact, if $x\in \Sigma'\backslash \Sigma$, then $\mathrm{dist}(x,\partial\Sigma)\leq \tilde{d}$, therefore
$x\in \overline{B_{\tilde{d}}(\partial\Sigma)}$. That is 
$\mathrm{dist}(x,\partial\Sigma')\leq \tilde{d}+\hat{d}$.
Finally, there exists a constant $C_3$, depending on the class $\mathcal{A}$ only, such that for any $t$, $0<t\leq 1$, we have
\begin{equation}\label{boundary33}
|\overline{B_{t}(\partial\Sigma)}|\leq C_3 t.
\end{equation}

We consider the following direct scattering problem.
Fixed $\Sigma\in\mathcal{A}$, for a fixed wavenumber $k>0$ and a fixed direction of propagation $v\in\mathbb{S}^{N-1}$, let the incident field $u^i$ be the corresponding time harmonic plane wave, that is $u^i(x)=\rme^{\rmi kx\cdot v}$, $x\in\mathbb{R}^N$. The incident wave is perturbed by the presence of the scatterer $\Sigma$ through a scattered wave, characterized by its corresponding scattered field $u^s$. The total field $u$ is the  
solution to the following exterior boundary value problem
 \begin{equation}\label{uscateq}
\left\{\begin{array}{ll}
\Delta u + k^2u=0 & \text{in }\mathbb{R}^N\backslash \Sigma\\
u=u^i+u^s & \text{in }\mathbb{R}^N\backslash \Sigma\\
B.C. & \text{on }\partial \Sigma\\
\displaystyle{\lim_{r\to \infty}r^{(N-1)/2}\left(\frac{\partial u^s}{\partial r}-\rmi ku^s\right)=0} & r=\|x\|,
\end{array}\right.
\end{equation}
where the last limit is the \emph{Sommerfeld radiation condition} and corresponds to the fact that
the scattered wave is radiating.
The boundary condition on the boundary of $\Sigma$ depends on the character  of the scatterer $\Sigma$. For instance, if $\Sigma$ is \emph{sound-soft}, then $u$ satisfies the following homogeneous Dirichlet condition
\begin{equation}\label{soundsoft}
u=0\quad\text{on }\partial \Sigma,
\end{equation}
whereas if $\Sigma$ is 
\emph{sound-hard} we have 
\begin{equation}\label{soundhard}
\nabla u\cdot \nu=0\quad\text{on }\partial \Sigma,
\end{equation}
that is a homogeneous Neumann condition.
Other conditions such as the impedance boundary condition or transmission conditions for penetrable scatterers may occur in the applications.

We recall that 
the Sommerfeld radiation condition holds uniformly with respect to all directions $\hat{x}=x/\|x\|\in\mathbb{S}^{N-1}$ and it implies that
the scattered field has the asymptotic behavior of an outgoing spherical wave, namely
\begin{equation}\label{asympt0}
u^s(x)=\frac{\rme^{\rmi k\|x\|}}{\|x\|^{(N-1)/2}}\left\{u_{\infty}(\hat{x})
+O\left(\frac{1}{\|x\|}\right)\right\},
\end{equation}
where $\hat{x}=x/\|x\|\in \mathbb{S}^{N-1}$ and
$u_{\infty}$ is the so-called \emph{far-field pattern} of $u^s$. In particular, the scattered field satisfies the following decay property for some positive constants $E$ and $R_1$
\begin{equation}\label{decay0}
|u^s(x)|\leq E\|x\|^{-(N-1)/2}\quad\text{for any }x\in\mathbb{R}^N\text{ so that }\|x\|\geq R_1.
\end{equation}

We refer to \cite{Wil} for further details, such as existence and uniqueness of the solution, on the direct scattering problem
\eqref{uscateq}. For an introduction to the corresponding inverse problems see for instance \cite{Col-Kre98,Isak06}.

Let us fix constants $0<\underline{k}<\overline{k}$ and let us denote, for any $N\geq 2$,
\begin{equation}\label{I_Ndefin}
I_N=\left\{\begin{array}{ll}
[\underline{k},\overline{k}]&\text{if }N=2,\\
(0,\overline{k}]&\text{if }N\geq 3.
\end{array}\right.
\end{equation}

\begin{prop}\label{l2bound}
Let us fix constants $0<\underline{k}<\overline{k}$ and let $I_N$ be defined as in \eqref{I_Ndefin}.
Let $\mathcal{A}$ be as defined at the beginning of the subsection.

Fixed $\Sigma\in\mathcal{A}$, $k\in I_N$, and $v\in\mathbb{S}^{N-1}$, let 
$u^i(x)=\rme^{\rmi kx\cdot v}$, $x\in\mathbb{R}^N$, and $u_{\Sigma,k,v}$ be the solution to \eqref{uscateq}, with boundary condition \eqref{soundsoft} or \eqref{soundhard},
and $u^s_{\Sigma,k,v}$ be its corresponding scattered field.

Then there exists a constant $E$, depending on $\mathcal{A}$ and
$I_N$ only, such that
\begin{equation}\label{uniformboundinfty}
|u_{\Sigma,k,v}(x)|\leq E\quad\text{for any }x\in \mathbb{R}^N\backslash\Sigma.
\end{equation}

Furthermore,
there exists a constant $E_1$, depending on the constant  $E$ in \eqref{uniformboundinfty}, $I_N$, $R$ and $N$ only, 
such that for any $\Sigma\in \mathcal{A}$, any $k\in I_N$, and any $v\in\mathbb{S}^{N-1}$
we have
\begin{equation}\label{udecayestimate}
|u^s_{\Sigma,k,v}(x)|\leq E_1\|x\|^{-(N-1)/2}\quad\text{for any }x\in\mathbb{R}^N\text{ so that }\|x\|\geq R+2.
\end{equation}
and
\begin{equation}\label{gradientdecayestimate}
\|\nabla u^s_{\Sigma,k,v}(x)\|\leq E_1\|x\|^{-(N-1)/2}\quad\text{for any }x\in\mathbb{R}^N\text{ so that }\|x\|\geq R+2.
\end{equation}
\end{prop}

\proof{.} First of all, we show that there exists a constant $E_0$, depending on $\mathcal{A}$ and
$I_N$ only, such that
\begin{equation}\label{uniformbound}
\|u_{\Sigma,k,v}\|_{L^2(B_{R+3}\backslash \Sigma)}\leq E\quad\text{for any }\Sigma\in \mathcal{A},\text{ any }k\in I_N,\text{ and any }v\in\mathbb{S}^{N-1}.
\end{equation}
This is an immediate consequence of Proposition~3.2 and Theorem~3.9 in \cite{Men-Ron} for the sound-hard case and of Lemma~3.5 in \cite{Ron03} for the sound-soft case.
Already from this first bound we can easily infer that \eqref{udecayestimate} and \eqref{gradientdecayestimate} hold true.

The main idea of the proof needed to 
improve the uniform $L^2$ bound in \eqref{uniformbound} to the uniform $L^{\infty}$ one contained in \eqref{uniformboundinfty} is the following.

Let $x\in \partial\Sigma$ and let us 
exploit Condition~\ref{pcondition}) of Definition~\ref{scatdefin}.
%, in its version  contained in Remark~\ref{suffcondremark}.
%
%assume that for some positive $\tilde{r}$, $\tilde{r}_1$ and $\tilde{C}$,
%and any connected component $U$ of
%$G\cap B_{\tilde{r}}(x)$, we can find an open set $U_1$, such that
%$U\cap B_{\tilde{r}_1}(x)\subset U_1\subset U$, and a bijective map $T:U_1\to (0,1)^N$ such that the following
%properties hold. The maps $T$ and $T^{-1}$ are locally Lipschitz and
%$\|D T\|$ and $\|DT^{-1}\|$ are bounded by $\tilde{C}$ almost everywhere. Furthermore, if we set
%$\Gamma=[0,1]^{N-1}\times\{1\}$, we require that $T^{-1}(\Gamma)=\partial U_1\cap \partial G$ and
%$T^{-1}(y)\in G$ for any $y\in [0,1]^N\backslash\Gamma$.
%
%Then 
By a change of variables, a reflection argument and standard regularity estimates, we infer that we can bound $|u|$ almost everywhere in $B_{r_2}(x)$ by a constant $\tilde{C}_1$,
where $r_2$ and $\tilde{C}_1$ depend on $r$, $r_1$, $\tilde{C}$ and the $L^2$ norm of $u$ which is bounded by \eqref{uniformbound}.

This procedure allows to estimate $|u|$ in a neighborhood of $\partial\Sigma$.
%, where $\tilde{K}$ is as in the proof of Proposition~\ref{suffcond}. Then we consider a point $x$ belonging to the boundary of some $K^i$, $i\in\{1,\ldots,M\}$, which are outside such a neighborhood. Then for some $r_1\leq r$ we have that $B_{r_1}(x)\cap \Sigma=B_{r_1}(x)\cap K^i$ and, by the assumptions on Lipschitz hypersurfaces, we have a bi-Lipschitz map $\Phi_x$ mapping $B_{r_1}(x)\backslash \Sigma$ onto
%$\Phi_x(B_{r_1}(x))\backslash \pi^+$
%where $\pi^+=\{y\in\mathbb{R}^N:\ y_N=0,\ y_{N-1}\geq 0\}$. Then it is easy to transform, in the required way, $\Phi_x(B_{r_1}(x))\backslash \pi^+$ onto a portion of a half-space and obtain the 
%required property. Thus we can bound $|u|$ in a neighborhood of the union of the boundaries of $K^i$, $i=1,\ldots,M$.
%Then we pass to points belonging to some $K^i$, $i\in\{1,\ldots,M\}$, which are outside such a neighborhood. With the same construction we can therefore bound $|u|$ in a whole neighborhood of $\partial \Sigma$.
Away from $\partial \Sigma$ the estimate is completely standard.\cvd

\bigskip

Let us fix $\Sigma$ and $\Sigma'$ belonging $\mathcal{A}$, $\mathcal{A}$ 
as defined at the beginning of the subsection. We also fix $k>0$ and a direction of propagation $v\in\mathbb{S}^{N-1}$.
Let $u$ be the solution to \eqref{uscateq} with boundary condition \eqref{soundsoft} or \eqref{soundhard}.
We denote by $u^s$ the corresponding scattered field and by $u_{\infty}$ its far-field pattern.
Moreover, $u'$, $(u^s)'$ and $u'_{\infty}$ denotes the same functions when $\Sigma$ is replaced by $\Sigma'$.
Finally, we fix positive $R_1$ and $\tilde{\rho}$ such that
$R+1+\tilde{\rho}\leq R_1$.

By Proposition~\ref{l2bound}, we have that
\begin{equation}\label{globalbound}
|u(x)|+|u'(x)|\leq E\quad\text{for any }x\in \mathbb{R}^3,
\end{equation}
where $E$ depends on $k$ and $\mathcal{A}$ only and it may be assumed to be greater than or equal to $1$, and $u$ and $u'$ are extended to $0$ on $\Sigma$ and $\Sigma'$, respectively.

Let us fix a point $x_0\in\mathbb{R}^N$ such that
$R+1+\tilde{\rho}\leq \|x_0\|\leq R_1$.
For a fixed $\varepsilon$, $0<\varepsilon\leq E$, let
\begin{equation}\label{errornear}
\|u-u'\|_{L^{\infty}(B_{\tilde{\rho}}(x_0))}\leq \varepsilon.
\end{equation}
We call $\varepsilon$ the \emph{near-field error with limited aperture}.

Then, let $\varepsilon_1>0$ be such that 
\begin{equation}\label{errornear2}
\|u-u'\|_{L^{\infty}(B_{\|x_0\|+\tilde{\rho}}\backslash \overline{B_{\|x_0\|-\tilde{\rho}})}}\leq \varepsilon_1.
\end{equation}
We call $\varepsilon_1$ the \emph{near-field error}.

Finally,
if
\begin{equation}\label{errorfar}
\|u_{\infty}-u'_{\infty}\|_{L^2(\mathbb{S}^{N-1})}\leq \varepsilon_0,
\end{equation}
$\varepsilon_0$ will be referred to as the \emph{far-field error}.

We investigate the relations among these errors. First of all, 
let us recall that a three-spheres inequality holds for the Helmholtz equation, provided the larger ball has a radius bounded by a constant $\tilde{\rho}_0$, $\tilde{\rho}_0$ depending on $k$ only,
see for instance \cite{Bru} or for a version suited to our case \cite[Lemma~3.5]{Ron08}
which we state here for the convenience of the reader.

\begin{lem}\label{3sphereslemma}
There exist positive  constants $\tilde{\rho}_0$, $C$ and $c_1$, $0<c_1<1$,
depending on $k$ only, such that for every 
$0<\rho_1<\rho<\rho_2\leq \tilde{\rho}_0$
and any function $u$ such that
$$\Delta u+k^2u=0\quad\text{in }B_{\rho_2},$$
we have, for any $s$, $\rho<s<\rho_2$,
\begin{equation}\label{3spheres}
\|u\|_{L^{\infty}(B_{\rho})}\leq C(1-(\rho/s))^{-N/2}
\|u\|^{1-\beta}_{L^{\infty}(B_{\rho_2})}\|u\|^{\beta}_{L^{\infty}(B_{\rho_1})},
\end{equation}
for some $\beta$ such that
\begin{equation}\label{3spherescoeff}
c_1\left(\log(\rho_2/s)\right)\big/\left(\log(\rho_2/\rho_1)\right)
\leq\beta\leq 1-c_1\left(\log(s/\rho_1)\right)\big/\left(\log(\rho_2/\rho_1)\right).
\end{equation}
\end{lem}

By an iterated application of this three-spheres inequality,
we have that there exist positive constants $C$ and $\alpha$, $0<\alpha<1$, depending on $E$, $\tilde{\rho}$, $R_1$ and $k$ only, such that
\begin{equation}\label{ballannulus}
\varepsilon\leq \varepsilon_1\leq C\varepsilon^{\alpha}.
\end{equation}

Moreover, 
there exist positive constants $\tilde{\varepsilon}_0< 1/\rme$ and $C_1$, depending on 
$E$, $R$, $\tilde{\rho}$, $R_1$ and $k$ only,
such that if $0<\varepsilon_0\leq\tilde{\varepsilon}_0$ then
\begin{equation}\label{fartonear}
\|u-u'\|_{L^{\infty}(B_{\|x_0\|+\tilde{\rho}}\backslash \overline{B_{\|x_0\|-\tilde{\rho}})}}\leq
\eta_1(\varepsilon_0)=\exp\left(-C_1(-\log\varepsilon_0)^{1/2}\right)
\end{equation}
that is
\begin{equation}\label{fartonear2}
\varepsilon\leq \varepsilon_1\leq \eta_1(\varepsilon_0)=\exp\left(-C_1(-\log\varepsilon_0)^{1/2}\right).
\end{equation}
This estimate follows immediately by the results in \cite{Isak92} for $N=3$ and with an easy modification for any other $N\geq 2$, see for instance Theorem~4.1 in \cite{Ron-Sin}.
%In this case, we may also consider $R_1$ and $\tilde{\rho}$ as depending on the other a priori data.

If we wish to reduce to obstacles only, we use the class of admissible obstacles $\mathcal{A}_{obst}$. In particular, we set $\mathcal{A}_{obst}=\hat{\mathcal{D}}_{obst}(r,L,R)$.

It is important to notice that $\delta$ in this case may be chosen to be as in \eqref{unifconnLipschitz}, therefore $\delta^{-1}$ may be chosen to be $C_2d$ for any $d$, for some constant $C_2$ depending on $r$, $L$ and $R$ only, that is \eqref{distancesopposite} becomes
\begin{equation}\label{distancesoppositebis}
C_1 d\leq C_1\hat{d}\leq \tilde{d}\leq C_2d\leq C_2\hat{d}.
\end{equation}

Finally, if we wish to use classes of admissible polyhedral scatterers or obstacles, fixed the size parameter $h>0$, we use
$\mathcal{A}^h=\hat{\mathcal{B}}_{scat}^h(r,L,R,r_1,\tilde{C},\omega,\delta)$ for 
general scatterers and
$\mathcal{A}_{obst}^h=\hat{\mathcal{D}}_{obst}^h(r,L,R)$ for obstacles.

\section{The main stability results}\label{resultssec}

In this section we present our stability results for the determination of sound-hard polyhedral scatterers. We distinguish them with respect to the number of scattering measurements used.

In this section we fix positive constants $r$, $L$ and $R$, $0<r_1< r$ and $\tilde{C}>0$.
Let us also fix $\omega:(0,+\infty)\to (0,+\infty)$ and $\delta:(0,+\infty)\to (0,+\infty)$ two  nondecreasing left-continuous functions.
We recall that $\omega(t)\leq t$, that
$\lim_{t\to+\infty}\omega(t)$ is equal to a finite real number which we call $\omega(+\infty)$,
and that $\delta(t)\leq t$ for any $t>0$. We fix the wavenumber $k>0$.
Finally, we fix positive $R_1$ and $\tilde{\rho}$ such that
$R+1+\tilde{\rho}\leq R_1$.
We refer to these constants and functions, including $N$, as the \emph{a priori data} and we let $\mathcal{A}=\hat{\mathcal{B}}_{scat}(r,L,R,r_1,\tilde{C},\omega,\delta)$ be the class of scatterers defined in Definition~\ref{scatdefin1}.
As before, for any fixed $h>0$ we call $\mathcal{A}^h=\hat{\mathcal{B}}_{scat}^h(r,L,R,r_1,\tilde{C},\omega,\delta)$ the set of scatterers $\Sigma\in \mathcal{A}$ such that $\Sigma$ is polyhedral with constants $h$ and $L$.

We call $\eta:(0,1/\rme)\to (0,+\infty)$ the following function
\begin{equation}\label{contmodul}
\eta(s)=\exp(-(\log(-\log s))^{1/2})\quad\text{for any }s,\ 0<s< 1/\rme.
\end{equation}

\subsection{Polyhedral scatterers with $N$ measurements}\label{Nmeasstatement}

 We fix $N$ linearly independent unit vectors $v_1,\ldots,v_N$.
 We notice that, given $N$ linearly independent unit vectors $v_1,\ldots,v_N$,
there exists a positive constant $a_0$, depending on the vectors $v_1,\ldots,v_N$, 
such that
\begin{equation}\label{independency}
\min_{ \nu \in \mathbb{S}^{N-1}} \left\{\max_{j \in \{1,\ldots,N\}} |v_j \cdot \nu|\right\} \geq a_0.
\end{equation}
In fact, $\max_{j \in \{1,\ldots,N\}} |v_j \cdot \nu|$ is a continuous function of $\nu\in \mathbb{S}^{N-1}$
which never vanishes.

We also fix a point $x_0\in\mathbb{R}^N$ such that
$R+1+\tilde{\rho}\leq \|x_0\|\leq R_1$.

\begin{teo}\label{mainteoN}
Let $N\geq 2$.
Fix $h>0$.
Let $\Sigma$, $\Sigma'$ belong to $\mathcal{A}^h$ and let $d$ be defined as in \eqref{ddefin}.
For any $j=1,\ldots,N$, let 
$u^i(x)=\rme^{\rmi kx\cdot v_j}$, $x\in\mathbb{R}^N$, and let $u_j$ be the solution to \eqref{uscateq} with boundary condition \eqref{soundhard} and $u'_j$ be the solution to the same problem with $\Sigma$ replaced by $\Sigma'$.

%There exists a constant $\hat{\varepsilon}>0$, depending on the a priori data only,
If
\begin{equation}\label{error2}
\max_{j=1,\ldots,N}\|u_j-u'_j\|_{L^{\infty}(B_{\tilde{\rho}}(x_0))}\leq\varepsilon
\end{equation}
for some $\varepsilon\leq 1/(2\rme)$,
then for some positive constant $C$ depending on the a priori data and on $a_0$ only, and not on $h$, we have
\begin{equation}\label{penultima}
\min\{d,h\}\leq 2\rme R(\eta(\varepsilon))^C.
\end{equation}
Therefore,
\begin{equation}\label{ultima1}
d\leq 2\rme R(\eta(\varepsilon))^C,
\end{equation}
provided $\varepsilon\leq \hat{\varepsilon}(h)$ where
\begin{equation}\label{epsilon0defin}
\hat{\varepsilon}(h)=
\min\bigg\{1/(2\rme),\eta^{-1}\bigg(\Big(\frac{h}{2\rme R}\Big)^{1/C}\bigg)\bigg\}.
\end{equation}
\end{teo}

With little modification, we obtain exactly the same stability result if $\Sigma$ and $\Sigma'$ are sound-soft scatterers instead of sound-hard ones, even if we reduce the number of measurements from $N$ to $1$. That is, as a byproduct of this work, 
we can significantly extend Theorem~4.1 in \cite{Ron08} to a much more general class of scatterers, namely $\mathcal{A}^h$, and
to any dimension $N\geq 2$.
We state such result in the following theorem.

\begin{teo}\label{mainteosoft}
Let $N\geq 2$.
Fix $h>0$.
Let $\Sigma$, $\Sigma'$ belong to $\mathcal{A}^h$ and let $d$ be defined as in \eqref{ddefin}. Let us fix $v\in \mathbb{S}^{N-1}$ and let $u^i(x)=\rme^{\rmi kx\cdot v}$, $x\in\mathbb{R}^N$. Let $u$ be the solution to \eqref{uscateq} with boundary condition \eqref{soundsoft} and $u'$ be the solution to the same problem with $\Sigma$ replaced by $\Sigma'$.

%There exists a constant $\hat{\varepsilon}>0$, depending on the a priori data only,
%such that i
If
\begin{equation}\label{error2s}
\|u-u'\|_{L^{\infty}(B_{\tilde{\rho}}(x_0))}\leq\varepsilon
\end{equation}
for some $\varepsilon\leq 1/(2\rme)$,
then for some positive constant $C$ depending on the a priori data only, and not on $h$, we have
\begin{equation}\label{penultimas}
\min\{d,h\}\leq 2\rme R(\eta(\varepsilon))^C.
\end{equation}
\end{teo}

\subsection{Polyhedral obstacles with fewer measurements}\label{lessmeasstatement}

It is well-known that in general $N-1$ scattering measurements may not be enough
to uniquely determine a polyhedral sound-hard screen. However, if we limit ourselves to polyhedral obstacles, that is to polyhedra, then a single measurement is enough, see \cite{Els-Yam1,Els-Yam2}.

Here we restrict ourselves to obstacles and we aim to obtain corresponding stability estimates with a minimal number of scattering measurements.

For technical reasons we limit ourselves to $N=2$ or $N=3$ only. Let us then
fix $N\in \{2,3\}$ and
positive constants $r$, $L$ and $R$. We fix the wavenumber $k>0$.
Finally, we fix positive $R_1$ and $\tilde{\rho}$ and a point $x_0\in\mathbb{R}^N$ such that
$R+1+\tilde{\rho}\leq \|x_0\|\leq R_1$.

We let $\mathcal{A}_{obst}=\hat{\mathcal{D}}_{obst}(r,L,R)$
be the class of scatterers defined in Definition~\ref{scatdefin1}.
For any fixed  $h>0$, we call $\mathcal{A}_{obst}^h=\hat{\mathcal{D}}_{obst}^h(r,L,R)$ the set of obstacles $\Sigma\in \mathcal{A}_{obst}$ such that $\Sigma$ is polyhedral with constants $h$ and $L$.

We recall that, for some constants and functions depending on $r$, $L$ and $R$ only, we have
$\hat{\mathcal{D}}_{obst}(r,L,R)\subset \hat{\mathcal{B}}_{scat}(\tilde{r},\tilde{L},R,r_1,\tilde{C},\omega,\delta)$.
Therefore, in this case we may set the constants
$r$, $L$, $R$, $k$, $R_1$ and $\tilde{\rho}$, including $N$,
as the \emph{a priori data}.

%It is important to notice that $\delta$ in this case may be chosen to be as in \eqref{unifconnLipschitz}, therefore $\delta^{-1}$ may be chosen to be $C_2d$ for any $d$, for some constant $C_2$ depending on $r$, $L$ and $R$ only, that is \eqref{distancesopposite} becomes
%\begin{equation}\label{distancesoppositebis}
%C_1 d\leq C_1\hat{d}\leq \tilde{d}\leq C_2d\leq C_2\hat{d}.
%\end{equation}

We begin by investigating an intermediate case, namely the one with $N-1$ scattering measurements. %In the next section we shall tackle the single measurement case.

\begin{teo}\label{mainteoN-1}
Let $N=2,3$.
Fix $h>0$. 
Let $\Sigma$, $\Sigma'$ belong to $\mathcal{A}_{obst}^h$ and let $d$ be defined as in \eqref{ddefin}.

If $N=2$, let us fix $v\in \mathbb{S}^{1}$ and let $u^i(x)=\rme^{\rmi kx\cdot v}$, $x\in\mathbb{R}^2$. Let $u$ be the solution to \eqref{uscateq} with boundary condition \eqref{soundhard} and $u'$ be the solution to the same problem with $\Sigma$ replaced by $\Sigma'$. We let $\varepsilon>0$ be such that
\begin{equation}\label{error3}
\|u-u'\|_{L^{\infty}(B_{\tilde{\rho}}(x_0))}\leq \varepsilon.
\end{equation}

If $N=3$, let us fix $v_1$, $v_2\in  \mathbb{S}^{2}$, with $|v_1\cdot v_2|=b_0<1$.
For any $j=1,2$, let 
$u^i(x)=\rme^{\rmi kx\cdot v_j}$, $x\in\mathbb{R}^3$, and let $u_j$ be the solution to \eqref{uscateq} with boundary condition \eqref{soundhard} and $u'_j$ be the solution to the same problem with $\Sigma$ replaced by $\Sigma'$.
We let $\varepsilon>0$ be such that
\begin{equation}\label{error4}
\max_{j=1,2}\|u_j-u'_j\|_{L^{\infty}(B_{\tilde{\rho}}(x_0))}\leq\varepsilon.
\end{equation}

There exists a constant $\hat{\varepsilon}_1(h)$, $0<\hat{\varepsilon}_1(h)
\leq 1/(2\rme)$, depending on the a priori data, on $b_0$ if $N=3$, and on $h$
only,
such that if $\varepsilon\leq \hat{\varepsilon}_1(h)$,
then for some positive constants $A_1$, depending on the a priori data only, and $C$, depending on the a priori data, on $b_0$ if $N=3$, and on $h$ only, we have
\begin{equation}\label{penultima1}
d\leq A_1(\eta(\varepsilon))^C.
\end{equation}
\end{teo}

%By the arguments developed at the beginning of the previous section and in \eqref{relationdistances}, it is an elementary exercise to infer the corresponding estimates of Theorem~\ref{mainteo1} if we replace $d$ with 
%$d_H(\Sigma,\Sigma')$ or $d_H(\partial\Sigma,\partial\Sigma')$ or the near-field error $\varepsilon_1$ with 
%a near-field error with limited aperture $\varepsilon$ or a far-field error $\varepsilon_0$ on the corresponding solutions.
%
%\subsection{The single measurement case}\label{singlestatement}

The main difference with respect to the sound-hard case with $N$ measurements
or the sound-soft case is that here we do not have an explicit dependence of 
$\hat{\varepsilon}_1(h)$ from $h$, which in Theorems~\ref{mainteoN} and \ref{mainteosoft} is given by
\eqref{epsilon0defin}, and that the constant $C$ depends, again in a rather implicit way, on $h$ too.

We finally consider the case of a single scattering measurement.
We restrict here to $N=3$, since $N=2$ is clearly covered by the previous theorem.

\begin{teo}\label{mainteo1}
Let $N= 3$.
Fix $h>0$.
Let $\Sigma$, $\Sigma'$ belong to $\mathcal{A}^h_{obst}$ and let $d$ be defined as in \eqref{ddefin}. Let us fix $v\in \mathbb{S}^2$ and let $u^i(x)=\rme^{\rmi kx\cdot v}$, $x\in\mathbb{R}^3$. Let $u$ be the solution to \eqref{uscateq} with boundary condition \eqref{soundhard} and $u'$ be the solution to the same problem with $\Sigma$ replaced by $\Sigma'$.

There exists a constant $\hat{\varepsilon}_2(h)$,
$0<\hat{\varepsilon}_2(h)\leq\hat{\varepsilon}_1(h)
\leq 1/(2\rme)$,
 depending on the a priori data
and on $h$ only,
such that if
\begin{equation}\label{errorN1}
\|u-u'\|_{L^{\infty}(B_{\tilde{\rho}}(x_0))}\leq\varepsilon
\end{equation}
for some $\varepsilon\leq\hat{\varepsilon}_2(h)$,
then for some positive constants $A_2\geq A_1$, depending on the a priori data only, and $C_1\leq C$, depending on the a priori data and on $h$ only, we have
\begin{equation}\label{penultima111}
d\leq A_2(\eta(\varepsilon))^{C_1}.
\end{equation}
\end{teo}

\begin{oss}\label{diffmeasrem}
We finally notice that, by the arguments developed in the previous section, we can easily infer corresponding estimates of Theorems~\ref{mainteoN} and \ref{mainteosoft},
and of Theorems~\ref{mainteoN-1} and \ref{mainteo1},
if we replace $d$ with 
$\tilde{d}=d_H(\Sigma,\Sigma')$ or $\hat{d}=d_H(\partial\Sigma,\partial\Sigma')$ or the near-field error with limited aperture $\varepsilon$ with 
a near-field error $\varepsilon_1$ or a far-field error $\varepsilon_0$ on the corresponding solutions.
In the first case, it is just enough to use \eqref{distancesopposite}, with $\delta^{-1}$ defined as in \eqref{delta-1def} and $C_1>0$ depending on the a priori data only,  for the first two theorems, and to use \eqref{distancesoppositebis}, with $C_1>0$ and $C_2$ depending on the a priori data only, for the second two theorems. For the second case, by \eqref{fartonear2}, we have exactly the same results if we replace $\varepsilon$ with the near-field error $\varepsilon_1$. If we use the far-field error $\varepsilon_0$ instead, then we need to replace $\varepsilon$ with $\eta_1(\varepsilon_0)$, $\eta_1$ as in
\eqref{fartonear2}, noting that in this case we can choose $\tilde{\rho}$ and $R_1$ as depending on the other a priori data.
\end{oss}

\section{The general geometric construction}\label{geomconstrsec}

In this section we assume that the assumptions of Theorem~\ref{mainteoN} are satisfied. The a priori data will be the one defined at the beginning of Section~\ref{resultssec}.

Moreover, for the whole section we shall fix $j\in\{1,\ldots,N\}$ and we shall consider the solutions with respect to the incident direction of propagation $v=v_j$, therefore the subscript $j$ will be always dropped.
 
We call $H$ the connected component of $G\cap G'$, where $G'=\mathbb{R}^N\backslash \Sigma'$,
such that $\mathbb{R}^N\backslash \overline{B_R}$
is contained in $H$.
We shall also use the following definition.

\begin{defin}\label{regchaindefin}
A sequence of balls $B_{\rho_i}(z_i)$, $i=0,\ldots,n$, forms a
\emph{regular chain}, with respect to an open set $G$, with constants $0<a_1<a_2<a_3<1<a_4$ 
if the following properties are satisfied
\begin{enumerate}[i)]
\item for any $i=0,1,\ldots,n$, $B_{a_4\rho_i}(z_i)\subset G$\textnormal{;}
\item for any $i=1,\ldots,n$, we have $\rho_i\leq \rho_{i-1}$ and
$B_{a_1\rho_{i}}(z_{i})\subset B_{a_2\rho_{i-1}}(z_{i-1})$
and, for any $i=0,\ldots,n-1$, we have
$B_{a_1\rho_{i}}(z_{i})\subset B_{a_3\rho_{i+1}}(z_{i+1})$.
\end{enumerate}
\end{defin}

We have the following lemmas with simple proofs that we leave to the reader.

\begin{lem}\label{regchainlemma}
Let $U_1$ and $\tilde{U}_1$ be two open sets and let $T:U_1\to \tilde{U}_1$ be a bi-$W^{1,\infty}$ mapping with constant $\tilde{C}$.

Let $B_{\tilde{\rho}_i}(\tilde{z}_i)$, $i=0,\ldots,n$, be a regular chain with respect to $\tilde{U}_1$ with constants $0<\tilde{a}_1<\tilde{a}_2<\tilde{a}_3<1<\tilde{a}_4$.
Then, if we call $z_i=T^{-1}(\tilde{z}_i)$, $\rho_i=\tilde{\rho}_i/\tilde{C}$, $i=0,\ldots,n$, and $a_1=\tilde{a}_1$, $a_4=\tilde{a}_4$, we have that $B_{\rho_i}(z_i)$, $i=0,\ldots,n$, 
is a regular chain with respect to $U_1$ with constants $0<a_1<a_2<a_3<1<a_4$ provided
$$a_1<\tilde{C}^2\tilde{a}_2\leq a_2<a_3<1\quad\text{and}\quad
a_1<\tilde{C}^2\tilde{a}_3\leq a_3<1.$$
\end{lem}

\begin{lem}\label{regchainlemma2}
Let $\mathcal{C}$ be an open cone with amplitude $\theta$, $0<\theta<\pi/2$, and radius $r$.
For simplicity we assume that its vertex is in the origin and that its bisecting
vector is $e_N$. We set $0<a_1<a_2<a_3<1<a_4$ and we call $0<c_1=\sin(\theta)/a_4<1$. Given $c_2$, $0<c_2\leq c_1$, we fix
$z_0=(r/2)e_N$ and $\rho_0$,
$c_2 (r/2)\leq \rho_0\leq c_1(r/2)$.

We can construct a regular chain $B_{\rho_i}(z_i)$, $i=0,\ldots,n$, 
with respect to $\mathcal{C}$ with constants $0<a_1<a_2<a_3<1<a_4$, in the following way.
For any $i=0,\ldots,n$, we can choose $z_n= b^n(r/2)  e_N$ and $\rho_n=  b^n \rho_0$
provided the constant $b$ satisfies 
\begin{equation}
0<\max\left\{\frac{1-a_2c_2}{1-a_1c_2},\frac{1+a_1c_2}{1+a_3c_2}\right\}\leq b<1.
\end{equation}
\end{lem}

We now proceed to describe the geometric construction needed for the proof of Theorem~\ref{mainteoN}, and of the other stability results as well.
We divide the construction into several steps, proving alongside their corresponding estimates.
Without loss of generality, up to swapping $\Sigma$ with $\Sigma'$,
we can find $x_1\in \partial\Sigma'\backslash\Sigma$ such that
$d=\mathrm{dist}(x_1,\partial\Sigma)=\mathrm{dist}(x_1,\Sigma)$. 

\subsubsection{Step I: from $x_0$ to $x_1$}

We construct a sequence of balls $B_{\rho_i}(z_i)$, $i=\ldots,-n,-(n-1),\ldots,0,\ldots,n_0$ forming a regular chain with respect to $G$,
with constants $0<a_1<a_2<a_3<1<a_4=8$ and $\rho_0$ depending on the a priori data only,
 and such that the following conditions are satisfied.

First, 
$z_0=x_0$ and $z_{n_0}=x_1$.
Second, $\rho_0$ is a positive constant, depending on the a priori data only,
such that
$16\rho_0\leq \min\{\tilde{\rho},\tilde{\rho}_0,r_1/\tilde{C}\}$, where $\tilde{\rho}_0$ is the positive constant depending on $k$ only that bounds the radius of balls where the three-spheres inequality of Lemma~\ref{3sphereslemma} holds.
On the other hand, 
 $\rho_{n_0}= s_0d$, where $s_0$ is a positive constant depending on the a priori data only.
Third, for any $n=1,2,\ldots$, we pick $z_{-n}=x_0+n(\rho_0/4)(x_0/\|x_0\|)$ and $\rho_{-n}=\rho_0$.
Finally $n_0$ is bounded by a constant, depending on the a priori data only, times
$\log(2\rme R/d)$.

The sequence is constructed as follows. Let $y_1$ be a point of $\partial\Sigma$ such that $|x_1-y_1|=d$. We recall that $B_d(x_1)\subset G$.

%Let $c$, $0<c<1/3$, be a constant depending on the a priori data only such that
%$\omega((1-c)r_1)\geq 1/2$.
If $d\geq r_1/3$, then we use the exterior connectedness property of $\Sigma$ and may easily construct such a chain keeping the radius $\rho_n=\rho_0$ for any $n\leq n_0$, that is simply constant and depending on the a priori data only. In this case we easily infer that $n_0$ is bounded by a constant depending on the a priori data only as well.

If instead $d\leq r_1/3$, we proceed in the following way. Let $U$ be the connected component of $G\cap B_{r_1}(y_1)$ containing $B_d(x_1)$.
In particular we have $y_1\in\partial U$. By Condition~\ref{pcondition}) of Definition~\ref{scatdefin} applied to $y_1$, we have the transformation $T:Q\to U'$
and we consider the point $\tilde{x}_1=T^{-1}(x_1)$.
We call $\tilde{y}_2$ the point in $Q$ such that $\tilde{y}_2=(\tilde{x}_1',3/4)$
and $y_2=T(\tilde{y}_2)$. We have that $B_{d/\tilde{C}}(\tilde{x}_1)\subset Q$ and that
$\mathrm{dist}(\tilde{x}_1,\partial Q\backslash \Gamma)\geq \omega(2r_1/3)$.
In particular,
$\|\tilde{x}'_1\|\leq 1-\omega(2r_1/3)$. We conclude that $\mathrm{dist}(y_2,\Sigma)$ is greater than or equal to a constant depending on the a priori data only.
By the exterior connectedness property of $\Sigma$ we construct such a chain
first from $x_0$ to $y_2$, keeping the radius constant and depending on the a priori data only. We notice that this part requires a number of balls that may be bounded by a constant depending on the a priori data only.
In order to proceed from $y_2$ to $x_1$, we use Lemma~\ref{regchainlemma2} to construct a regular chain in $Q$, with suitable constants, connecting $\tilde{y}_2$ to $\tilde{x}_1$. Then our chain in $G$ is obtained by using Lemma~\ref{regchainlemma} and we easily infer that
the number of elements of such a chain may be bounded by a constant, depending on the a priori data only, times $\log(2\rme R/d)$, therefore the claim is proved.

Let us finally notice that here the geometric construction is different from that of \cite{Ron08}. It is actually more general and more complicated and allows us to consider a wider class of admissible scatterers.

Starting from $z_0=x_0$, we take $j_1\in\{1,\ldots,n_0\}$
such that,
for any $i=0,1,\ldots,j-1$, $B_{\rho_i}(z_i)\subset H$ and
$B_{\rho_{j_1}}(z_{j_1})\cap \Sigma'\neq\emptyset$.
We apply the three-spheres inequality of Lemma~\ref{3sphereslemma} as follows. For any $i=0,1,\ldots,j-1$,
\begin{multline*}
\|u-u'\|_{L^{\infty}(B_{a_1\rho_{i+1}}(z_{i+1}))}\leq \|u-u'\|_{L^{\infty}(B_{a_2\rho_{i}}(z_{i}))}\leq\\
\leq C\|u-u'\|_{L^{\infty}(B_{\rho_{i}}(z_{i}))}^{1-\beta_i}\|u-u'\|_{L^{\infty}(B_{a_1\rho_{i}}(z_{i}))}^{\beta_i}.
\end{multline*}
where any $\beta_i$, $i=0,\ldots,j-1$, satisfies
$$0<a\leq\beta_i\leq b<1$$
with $a$ and $b$ depending on $k$ only.

If $\beta_i$, $i=0,1,2,\ldots$, are positive constants, we shall use the following notation
for any $j=0,1,2,\ldots$
$$\mathcal{B}_{j}=\sum_{r=0}^j\prod_{i=r}^j\beta_i,\quad\Gamma_j=\prod_{i=0}^j\beta_i.$$

Recalling that $\|u-u'\|_{L^{\infty}(\mathbb{R}^3)}\leq E$ and that $\|u-u'\|_{L^{\infty}(B_{\tilde{\rho}}(x_0))}\leq \varepsilon$, and by iterating the previous estimate, we obtain
\begin{equation}\label{firstest}
\|u-u'\|_{L^{\infty}(B_{a_1\rho_{j_1}}(z_{j_1}))}%=\varepsilon_2
\leq
 C^{1+\mathcal{B}_{j_1-1}}E^{1-\Gamma_{j_1-1}}\varepsilon^{\Gamma_{j_1-1}}.
\end{equation}

\subsubsection{Step II: towards the cell and back}

We call $\hat{h}=\min\{d,h\}$.
Following the notation of the previous step, we set $z=z_{j_1}$ and $\rho=\rho_{j_1}$.
Then, $B_{a_1\rho}(z)\subset H$, $B_{a_4\rho}(z)\subset G$, with $a_4=8$, and there exists
$w\in \partial\Sigma'$ such that $\|z-w\|<\rho$ and $B_{\|z-w\|}(z)\subset H$.
Let $U$ be
the connected component of $G'\cap B_{r_1}(w)$ containing $z$. Clearly $w\in\partial U$. Let $\mathcal{C}'$ be one of the cells of $\partial\Sigma'$ such that $w\in \mathcal{C}'$ and $\mathcal{C}'\cap B_{r_1/2}(w)\subset \partial U$.

We call $\Pi'$ the plane containing $\mathcal{C}'$ and, up to a rigid change of coordinates, without loss of generality, we assume
$\Pi'=\{y\in\mathbb{R}^N:\ y_N=0\}$.
By the properties of $\mathcal{C}'$, $\mathcal{C}'$ satisfies a uniform cone property, namely there exists $\omega_1\in \mathbb{S}^{N-1}\cap \Pi$ and constants $c_1$, $0<c_1\leq 1$, and $\theta$, $0<\theta<\pi/2$,
depending on $L$ and $R$ only, such that
$\mathcal{C}(w,\omega_1,c_1\hat{h},\theta)\cap\Pi'\subset \mathcal{C}'\cap B_{r_1/2}(w)\subset \partial U$.

By looking at the points on the bisecting line of $\mathcal{C}(w,\omega_1,c_1\hat{h},\theta)$, we may find
$w_1$ on this line, that is $w_1=w+s_1\hat{h}\omega_1$,
such that $B_{s_2\hat{h}}(w_1)\subset  B_{r_1/2}(w)\cap B_{7\rho/(4\tilde{C}^2)}(w)$,
and $B_{s_2\hat{h}}(w_1)\cap\Pi'\subset \mathcal{C}'\cap B_{r_1/2}(w)\subset \partial U$, for some positive $s_1$ and $s_2$ depending on the a priori data only.

We claim that
there exists $s_3$, $0<s_3\leq s_2$, depending on the a priori data only, such that,
up to changing the orientation of $e_N$, we have
$B^+_{s_3\hat{h}}(w_1)\subset U$.

The proof of this claim is the following. We apply Condition~\ref{pcondition}) of Definition~\ref{scatdefin} to $w$, and we consider the corresponding transformation $T:Q\to U'$.
For some $\varepsilon>0$, possibly taking $s_2$ slightly smaller but still depending on the a priori data only, we have that
$(\overline{B_{s_2\hat{h}}(w_1)}\cap\Pi')\times (0,\varepsilon]\subset U\cap B_{r_1/2}(w)$.
The function $T^{-1}$ restricted to such a set is bi-Lipschitz onto $A\subset Q$. We observe that $A$ has a positive distance, depending on $r_1$ and $\omega$ only, from $\partial Q\backslash \Gamma$.

It is not difficult to show that $T^{-1}$ can be extended to a function $\tilde{T}^{-1}$ in such a way that
$\tilde{T}^{-1}:(\overline{B_{s_2\hat{h}}(w_1)}\cap\Pi')\times [0,\varepsilon]\to \overline{A}$ is still bijective, with inverse $T$, and thus bi-Lipschitz. We conclude that also
$\tilde{T}^{-1}:(\overline{B_{s_2\hat{h}}(w_1)}\cap\Pi')\to \overline{A}\cap \Gamma$ is bijective, with inverse $T$, and thus bi-Lipschitz. Then let us fix $a$, $0<a\leq s_2\hat{h}$ such that
$B^+_a(w_1)\subset U\cap B_{r_1/2}(w)\subset U'$.
If $a$ can be chosen to be $s_2\hat{h}$, the claim is proved, otherwise we assume that
there exists $y\in \partial U'$ such that
$y\in B^+_{s_2\hat{h}}(w_1)$ and $\|y-w_1\|=a$.
Actually, $y\in \partial U\cap B_{r_1/2}(w)$.
Then we call $\tilde{w}_1=\tilde{T}^{-1}(w_1)$, thus $T(\tilde{w}_1)=w_1$, and, by a similar reasoning, extending $T^{-1}$ up to the closure of $B^+_a(w_1)$, we can find
$\tilde{y}\in \Gamma$ such that $T(\tilde{y})=y$
and $\|\tilde{y}-\tilde{w}_1\|\leq\tilde{C}a$. Moreover, we also have that
$\tilde{y}\not\in \overline{A}\cap \Gamma$.
But $B_{s_2\hat{h}/\tilde{C}}(\tilde{w}_1)\cap\Gamma\subset \overline{A}\cap \Gamma$, therefore $\tilde{C}a$ must be greater than $s_2\hat{h}/\tilde{C}$, that is,
$a>s_2\hat{h}/\tilde{C}^2$ and the claim is proved.

Then we notice that, through an even reflection, we can extend $u'$ on  $B_{s_3\hat{h}}(w_1)$ by setting $u'(y)=u'(T_{\Pi'}(y))$ for any $y\in B^-_{s_3\hat{h}}(w_1)$. In this way $u'$ solves the Helmholtz equation on the whole $B_{s_3\hat{h}}(w_1)$. We notice that $B_{s_3\hat{h}}(w_1)\subset G$, therefore on $H\cup
B_{s_3\hat{h}}(w_1)$ both $u$ and $u'$ are well defined and solve the Helmholtz equation.

We construct a regular chain, with respect to $H\cup
B_{s_3\hat{h}}(w_1)$ and with constants depending on the a priori data only,
satisfying the following properties. The first ball is centered in $z$ and it has radius less than or equal to $\rho$, whereas the last ball is $B_{s_5\hat{h}}(w_1)$ with $s_5$ depending on the a priori data only. Finally, the number of balls of such a chain is bounded by a constant depending on the a priori data only times $\log(2\rme R/\hat{h})$.

The argument is the following. By a reasoning similar to the one used before, we can find $\tilde{w}\in \Gamma$ such that $T(\tilde{w})=w$ and
$\|\tilde{w}_1-\tilde{w}\|\leq \tilde{C}s_1\hat{h}$.
Since $B_{7\rho}(w)\cap U\subset H$, we have $T(B_{7\rho/\tilde{C}}(\tilde{w})\cap Q)\subset H$. Since $s_1\geq s_3$, we have
$T^{-1}(B^+_{s_3\hat{h}}(w_1))\subset B_{2\tilde{C}s_1\hat{h}}(\tilde{w})$.
Without loss of generality, we require that $2\tilde{C}s_1\hat{h}\leq 7\rho/(2\tilde{C})$, that is, $T^{-1}(B^+_{s_3\hat{h}}(w_1))\subset B_{7\rho/(2\tilde{C})}(\tilde{w})\cap Q$.

Then we perform the following construction. Take $w_0$ on the segment connecting $z$ to $w$ such that $\|w_0-w\|= 7\rho/(4\tilde{C}^2)$, if this number is less than $\|z-w\|$.
Otherwise we pick $w_0=z$. 
We observe that $T^{-1}(w_0)\in B_{7\rho/(2\tilde{C})}(\tilde{w})\cap Q$.

We construct a regular chain of balls, with a number of balls bounded by a constant depending on the a priori data only, contained in $B_{\|z-w\|}(z)\subset H$ and connecting $z$ to $w_0$. We consider $w_2=w_1+s_4\hat{h}e_N$ so that
$B_{s_5\hat{h}}(w_2)\subset B^+_{s_3\hat{h}}(w_1)$. Then with a construction similar to one described before during Step I, which exploits the properties of the change of variables $T$, we can extend our regular chain, which is still contained in $H$, from $w_0$ till we connect to $w_2$. The number of steps required at this stage is of the order of a constant times $\log(2\rme R/\hat{h})$. Then we move along the segment connecting $w_2$ to $w_1$ and, with a finite number of steps depending only on the a priori data, we are able to reach $w_1$ and thus conclude the construction.

We notice again that, as in the first Step I, the construction developed here is much more general and much more involved than that used in \cite{Ron08}. Overcoming this technical difficulty is the key ingredient to obtain our results for a more general class of admissible scatterers than the one used in \cite{Ron08}.

Again by a repeated use of the three-spheres inequality, and by recalling
\eqref{firstest}, 
we obtain that
%\begin{equation}\label{secondest}
%\|u-u'\|_{L^{\infty}(B_{a_1s_3\hat{d}}(w_2))}=\varepsilon_3
%\leq C^{1+\mathcal{B}_{m-1}}(C_1d^{\alpha})^{1-\Gamma_{m-1}}\varepsilon_2^{\Gamma_{m-1}}
%\end{equation}
%where
%$$m\leq \tilde{C}\log(\rme d/\hat{d})$$
%and any $\beta_i$, $i=0,\ldots,m-1$, satisfies
%$0<a\leq\beta_i\leq b<1$, with
%$a$ and $b$ depending on $k$ only.
%
%Actually, coupling \eqref{secondest} with \eqref{firstest}, we also have
\begin{equation}\label{thirdest}
\|u-u'\|_{L^{\infty}(B_{a_1s_5\hat{h}}(w_1))}\leq
C^{1+\mathcal{B}_{n-1}}E^{1-\Gamma_{n-1}}\varepsilon^{\Gamma_{n-1}}
\end{equation}
where, for some constants $\tilde{C}_1$ and $a$, $b$, with $0<a<b<1$, 
depending on the a priori data only,
we have
\begin{equation}\label{nprop}
n\leq \tilde{C}_1\left(\log\left(\frac{2\rme R}{\hat{h}}\right)+\log\left(\frac{2\rme R}{d}\right)\right)\quad\text{and}\quad
a\leq\beta_i\leq b\text{ for any }i=0,\ldots,n-1.
\end{equation}

We then apply a reflection argument.
We call $\Pi_1=\Pi'$ the hyperplane containing the cell $\mathcal{C}'$. Moreover, $\nu_1$ will be the unit normal to $\Pi_1$ and 
$T_1= T_{\Pi_1}$ is the reflection in $\Pi_1$.
We define $\Sigma_1$ as the reflection of $\Sigma$ with respect to the plane $\Pi_1$, $G_1=\mathbb{R}^3\backslash \Sigma_1$, and $u_1$ as the even
reflection of $u$ with respect to the same plane $\Pi_1$,
namely for any $x\in\mathbb{R}^N$, we set $u_1(x)=u(T_{\Pi_1}(x))$.
Without loss of generality, and since $a_4=8$, we can assume that $B_{s_5\hat{h}}(w_1)\subset  B_{2\rho}(z)\subset B_{8\rho}(z)\subset G$, therefore
$B_{\rho}(z)\subset B_{3\rho}(w_1) \subset B_{4\rho}(w_1)\subset G\cap G_1$.
Both $u$ and $u_1$ satisfy the Helmholtz equation on $B_{4\rho}(w_1)$.
Notice that $\nabla u'\cdot \nu_1=0$ on
$\Pi_1\cap B_{s_5\hat{h}}(w_1)$, 
therefore $u'=u'\circ T_{\Pi_1}$ and
$$u-u_1=u-u'+u'-u_1=(u-u')-(u-u')\circ T_{\Pi_1}.$$
We can conclude, using \eqref{thirdest}, that
\begin{equation}
\|u-u_1\|_{L^{\infty}(B_{a_1s_5\hat{h}}(w_1))}\leq 2\|u-u'\|_{L^{\infty}(B_{a_1s_5\hat{h}}(w_1))}\leq
2C^{1+\mathcal{B}_{n-1}}E^{1-\Gamma_{n-1}}\varepsilon^{\Gamma_{n-1}}.
\end{equation}

Then by using the arguments of Step IV of Section~5 in \cite{Ron08}, we obtain 
that 
\begin{equation}\label{7est}
\|u-u_1\|_{L^{\infty}(B_{\rho_{j_1}}(z_{j_1}))}
\leq C^{1+\mathcal{B}_{n}}(2E)^{1-\Gamma_{n}}\varepsilon^{\Gamma_{n}}
\end{equation}
where $C\geq 1$ and $2E\geq 1$ are constants depending on the a priori data only,
\eqref{nprop} is satisfied and
$\beta_n$ satisfies
$$c_1\frac{\log(8/7)}{\log(c_2\rho_0/\hat{h})}
\leq\beta_n\leq 1-c_1+c_1\frac{\log(8/7)}{\log(c_2\rho_j/\hat{h})},$$
with $c_1$, $c_2$ depending on the a priori data only.

Finally, we call $w_1$ the first reflection point and $\Pi_1$ the first reflection hyperplane.

\subsubsection{Step III: returning back towards $x_0$ and infinity}
Let us beginning by fixing a constant $R_2\geq \max\{2R_1,4R\}$, depending on the a priori data and on $a_0$ only, such that 
$$E_1R_2^{-(N-1)/2}\leq ka_0/2$$
where $E_1$ is as in \eqref{gradientdecayestimate}.

Let us now consider the regular chain of balls $B_{\rho_i}(z_i)$, $i=\ldots,-n,\ldots,-1,0,1\ldots,j_1$, we have constructed in Step I.
We have that $B_{\rho_{j_1}}(z_{j_1})$ is contained in $G_1$.
We proceed backwards along the chain, until we find $j_2$, $ j_2<j_1$, such that, for any $i=j_2+1,\ldots,j_1$,
we have $B_{\rho_i}(z_i)\cap G_1=\emptyset$, whereas $B_{\rho_{j_2}}(z_{j_2})\cap G_1\neq\emptyset$. Then, we apply Step II
to $u$, $u_1$, $\Sigma$ and $\Sigma_1$. We find a second reflection point $w_2$ and a second reflection hyperplane $\Pi_2$, with unit normal $\nu_2$.
By reflection in such a hyperplane $\Pi_2$, from $\Sigma$ we obtain $\Sigma_2$ and from $u$ we obtain $u_2$. In a completely analogous way as in 
\eqref{7est}, we may estimate $\|u-u_2\|_{L^{\infty}(B_{\rho_{j_2}}(z_{j_2}))}$.

We repeat this procedure as many times as needed, until we reach $z_{-n_1}$, where
$n_1$ is an integer bounded by a constant depending on the a priori data and on $a_0$ only,  with $R_3=\|z_{-n_1}\|\geq 2R_2+2$. Fixed a hyperplane $\Pi$, to be decided later, that is passing through a point belonging to $\overline{B_{R_2+1}}$, and a point $\tilde{z}\in \partial B_{R_3}\cap \Pi$, by a regular chain of balls with constant radius $\rho_0$, we proceed from $z_{-n_1}$ along the boundary of $\partial B_{R_3}$
 towards the point $\tilde{z}\in \partial B_{R_3}\cap \Pi$.
 
 Before reaching $\tilde{z}$, we have done $M$ reflection procedures as in Step II, where $M$ is a positive integer bounded by $n_0$ plus a constant depending on the a priori data and on $a_0$ only.
 
We now distinguish between two cases. In the first, 
setting $\Sigma_0=\Sigma'$ and $u_0=u'$,
we assume that there exists a reflection point $w_n\in\partial\Sigma_{n-1}$, $1\leq n\leq M$, as above with $\|w_n\|\geq R_2+1$.
Then we have, since $\nabla u_{n-1}(w_n)\cdot \nu_n=0$,
$$\nabla u(w_n)\cdot\nu_n=\nabla (u-u_{n-1})(w_n)\cdot \nu_n.$$
Otherwise, in the second case, we can assume that the last reflection point $w_M$ is such that $\|w_M\|\leq R_2+1$ and, without loss of generality, we may pick $\Pi=\Pi_M$. Then, since $u_M=u\circ T_{\Pi_M}$,
we have $\nabla u_M(\tilde{z})\cdot\nu_M=-\nabla u(\tilde{z})\cdot\nu_M$, hence
$$2\nabla u(\tilde{z})\cdot\nu_M=\nabla (u-u_M)(\tilde{z})\cdot \nu_M.$$

In either cases, picking either $z=w_n$ or $z=\tilde{z}$,
we can prove the following lemma, see \cite[Section~5]{Ron08} for further details on the computations.

\begin{lem}\label{reflectionpointexistence}
We can find a point $z$, $\|z\|\geq R_2+1$ and a unit vector $\nu$
such that
\begin{equation}\label{quasiultimastima}
h|\nabla u(z)\cdot\nu|\leq C_0\varepsilon_2
\end{equation}
where 
for some $\beta_i$, $i=0,\ldots,n$,
\begin{equation}\label{finalest}
\varepsilon_2\leq C^{1+\mathcal{B}_{n}}(2E)^{1-\Gamma_{n}}\varepsilon^{\Gamma_{n}},
\end{equation}
with $C\geq 1$, $2E\geq 1$ and
$$n\leq \hat{C}\log(2\rme R/d)(\log(2\rme R/d)+\log(2\rme R/\hat{h})).$$
Furthermore, there are at most $M\leq \hat{C}_1\log(2\rme R/d)$ of these $\beta$ such that
$$\beta \geq c_1\frac{\log(8/7)}{\log(c_2\rho_0/\hat{h})}
%\leq \beta\leq 1-c_1+c_1\left(\log(8/7)\right)\big/\left(\log(c_2\rho_j/\hat{d})\right)$$
$$
and they are never consecutive ones,
and all the others satisfy
$0<a\leq\beta\leq b<1$.
Here $C_0$, $C$, $E$, $a$, $b$, $c_1$, $c_2$, $\hat{C}$ and $\hat{C}_1$ depend on the a priori data only.
\end{lem}

This lemma concludes the general geometric construction that is the basic step for proving our stability results.

\section{Proofs of the stability results}\label{proofssec}

In this section, using the geometric construction of the previous Section~\ref{geomconstrsec} as a starting point, we prove our stability results. For the $N$ measurements case the conclusion is straightforward, whereas if we consider less than $N$ measurements, we need to develop new arguments. 

\subsection{The $N$ measurements case}\label{Nmeassubs}

We conclude the proof of Theorem~\ref{mainteoN}, and thus also of Theorem~\ref{mainteosoft}.

\proof{ of Theorem~\ref{mainteoN}.}
By Lemma~\ref{reflectionpointexistence}, and using its notation, we have for any $j=1,\ldots,N$
$$\max_{j=1,\ldots,N}|\nabla u_j(z)\cdot\nu|\leq \frac{C_0}{h}\varepsilon_2.$$
Therefore, using Proposition~\ref{l2bound} and our choice of $R_2$, for any $j=1,\ldots,N$, we have
$$k|v_j\cdot \nu|-ka_0/2
%|\nabla u^i(z)\cdot\nu|-ka_0/2
\leq |\nabla u_j(z)\cdot\nu|\leq \frac{C_0}{h}\varepsilon_2.$$
Therefore, choosing one of the available incident waves we can infer that
\begin{equation}\label{almostfinal}
ka_0/2\leq \frac{C_0}{h}\varepsilon_2.
\end{equation}

Then, 
by straightforward although lengthy computations, see \cite[Section~5]{Ron08} for further details, the proof may be easily concluded.\cvd

\subsection{The $N-1$ measurements case}\label{N-1meassubs}

Here we assume that the hypotheses of Theorem~\ref{mainteoN-1} are satisfied.
Since we would like to keep our argument as general as possible, let us assume for the time being that $N\geq 2$ and that we have fixed $N-1$ linearly independent directions
$v_1,\ldots,v_{N-1}$.

We need the following lemma. 

\begin{lem}\label{obstacleslemma}
 There exists a constant $\tilde{a}_0>0$, depending on the a priori data only,
such that for any direction $v$, and any polyhedral $\Sigma\in \mathcal{A}_{obst}$, we can find a cell $\mathcal{C}$ in $\partial\Sigma$, with unit normal $\nu$, such that $|\nu\cdot v|\geq \tilde{a}_0$.
\end{lem}

\proof{.} Let us assume, by contradiction, that such a positive constant $\tilde{a}_0$ does not exist. Then we can find a sequence of polyhedral obstacles $\Sigma_n\in \mathcal{A}_{obst}$ and of directions $v_n$, $n\in\mathbb{N}$,
such that, for $\mathcal{H}^{N-1}$ almost any point $x$ of $\partial\Sigma_n$, we have $|\nu(x)\cdot v_n|\leq 1/n$. Without loss of generality, we can assume that, as $n\to\infty$, $\Sigma_n$ converges, in the Hausdorff distance, to $\Sigma\in\mathcal{A}_{obst}$
and that $v_n\to v\in \mathbb{S}^{N-1}$. We can conclude that for $\mathcal{H}^{N-1}$ almost any point $x$ of $\partial\Sigma$ we have $|\nu(x)\cdot v|=0$, which is impossible since $\Sigma$ is an obstacle.\cvd

\bigskip

We begin with the following interesting and not that difficult case. Let us consider the geometric construction of the previous section, in particular Lemma~\ref{reflectionpointexistence}. If the point $z$ defined there is a reflection point $w_n$, then a single measurement would be enough to obtain a stability result. In fact the following result holds.

\begin{prop}\label{easycaseprop}
Let $N\geq 2$.
Fix $h>0$. 
Let $\Sigma$, $\Sigma'$ belong to $\mathcal{A}_{obst}^h$ and let $d$ be defined as in \eqref{ddefin}.
Let us fix $v\in \mathbb{S}^{N-1}$ and let $u^i(x)=\rme^{\rmi kx\cdot v}$, $x\in\mathbb{R}^N$. Let $u$ be the solution to \eqref{uscateq} with boundary condition \eqref{soundhard} and $u'$ be the solution to the same problem with $\Sigma$ replaced by $\Sigma'$.
Let us assume that
$$
\|u-u'\|_{L^{\infty}(B_{\tilde{\rho}}(x_0))}\leq\varepsilon
$$
for some $\varepsilon\leq 1/(2\rme)$.

Let us further assume that the point $z$ defined in 
Lemma~\textnormal{\ref{reflectionpointexistence}} is a reflection point.
Then for some positive constant $C$ depending on the a priori data only, and not on $h$, we have
$$
\min\{d,h\}\leq 2\rme R(\eta(\varepsilon))^C.
$$
\end{prop}

\proof{.} 
The main idea of the proof is the following Let $z=w_n\in \Sigma_{n-1}$ be the reflection point such that $\|z\|\geq R_2+1$. We consider $\sigma$ the connected component of $\Sigma_{n-1}$ containing $z$ and we find, using Lemma~\ref{obstacleslemma},
a point $\tilde{z}\in \partial\sigma$ such that $|\nu(\tilde{z})\cdot v|\geq \tilde{a}_0$.

Clearly $\sigma$ is far away from $\Sigma$, therefore we are able to modify our regular chain by moving around $\partial\sigma$ till we get close to the point $\tilde{z}\in \partial\sigma$ and we can use such a point $\tilde{z}$ as a reflection point. Therefore the proof follows as in the proof of Theorem~\ref{mainteoN}, simply by replacing $a_0$ by $\tilde{a}_0$.\cvd

\bigskip

The difficult part arises when the assumptions of Proposition~\ref{easycaseprop} are not met, namely when in the geometric construction of the previous section we have $M$ reflection points, all of them contained in $B_{R_2+1}$. 

In the sequel, without loss of generality, we assume that in our geometric construction we have $M$ reflection points, all of them contained in $B_{R_2+1}$. Using the construction of the previous section, we can reach with a regular chain any point $\tilde{z}\in (
\overline{B_{2R_2+3}}\backslash B_{2R_2+2})
\cap \Pi_M$. Therefore, for any $j=1,\ldots,N-1$, we have
\begin{equation}\label{upperbound}
A_j=\max_{\tilde{z}\in (\overline{B_{2R_2+3}}\backslash B_{2R_2+2})
\cap \Pi_M }|\nabla u_j(\tilde{z})\cdot \nu|\leq C_0\varepsilon_2
\end{equation}
where $\nu=\nu_M$ is the unit normal to $\Pi_M$ and $C_0$ and $\varepsilon_2$ satisfy the same properties as those in Lemma~\ref{reflectionpointexistence}.

Let us illustrate what is the difficult point. In order to obtain our stability result we need to match the upper bound in \eqref{upperbound} with a corresponding lower bound.
Let us begin with the following remark. Let $v=v_j$, $j\in \{1,\ldots,N-1\}$, be one of the incident directions of propagation and, for the time being, let us drop the subscript $j$ from our solutions. Let us call
$$A=\max_{\tilde{z}\in (\overline{B_{2R_2+3}}\backslash B_{2R_2+2})
\cap \Pi_M }|\nabla u(\tilde{z})\cdot \nu|.$$
Can $A$ be equal to $0$? Indeed this can happen, although only in certain circumstances. Namely, we claim that $A=0$ if and only if 
$v\cdot \nu=0$ and
$\Sigma$ is symmetric with respect to the hyperplane $\Pi_M$.
One direction is obvious, let us show the more interesting one, that is  $A=0$ implies that $v\cdot \nu=0$ and
$\Sigma$ is symmetric with respect to the hyperplane $\Pi_M$.

In fact, if $A=0$, then $|\nabla u\cdot \nu|\equiv 0$ on 
$(\overline{B_{2R_2+3}}\backslash B_{2R_2+2})
\cap \Pi_M$ and, by unique continuation, we actually have that
$|\nabla u\cdot \nu|\equiv 0$ on $(\mathbb{R}^N\backslash \overline{B_R})\cap \Pi_M$. 
By the decay properties at infinity of $\nabla  u^s$, this may hold only if $v\cdot\nu=0$.
Moreover, we can easily infer that $u$ is even symmetric with respect to $\Pi_M$.
Let us call $\tilde{\Sigma}$ the complement of the unbounded connected component of $\mathbb{R}^{N}\backslash(\Sigma\cup T_{\Pi_M}(\Sigma))$. We have that $\tilde{\Sigma}$ is a polyhedral obstacle which is symmetric with respect to $\Pi_M$ and that $u$ solves \eqref{uscateq} with boundary condition \eqref{soundhard} also with $\Sigma$ replaced by $\tilde{\Sigma}$. 
By the uniqueness result for sound-hard polyhedral obstacles with a single measurement,
\cite{Els-Yam1,Els-Yam2}, we immediately infer that $\Sigma=\tilde{\Sigma}$ thus $\Sigma$ itself is symmetric with respect to $\Pi_M$.

Therefore, in order to bound $A$ away from $0$, we need to guarantee either that $\Sigma$ is not symmetric with respect to a hyperplane whose normal is orthogonal to $v$ or, if  $\Sigma$ is actually symmetric with respect to a hyperplane whose normal is orthogonal to $v$, that $\Pi_M$ is different from such a hyperplane.
As we shall see, if we use $N-1$ measurements, instead, in order to bound $\max_{j=1,\ldots,N-1}A_j$ away from zero, we need to guarantee either that $\Sigma$ is not symmetric with respect to a hyperplane whose normal is orthogonal to any $v_j$, $j=1,\ldots,N-1$, or, if  $\Sigma$ is actually symmetric with respect to such a hyperplane, that $\Pi_M$ is different from it.

Here lies the main difference between the $N-1$ measurements and $1$ measurement case. In fact, in the $N-1$ measurement case, for any obstacle $\Sigma$ there is at most one hyperplane  
whose normal is orthogonal to any $v_j$, $j=1,\ldots,N-1$, with respect to which
$\Sigma$ is symmetric. On the contrary, with only one measurements, for any obstacle $\Sigma$ there might be many hyperplanes 
whose normal is orthogonal to $v$ with respect to which
$\Sigma$ is symmetric. This is the main reason why the $N-1$ measurements case is
relatively simpler and the corresponding result is somewhat stronger. 

We then first consider the $N-1$ measurements case, leaving the $1$ measurement case to the next subsection.
Another difficulty is that we not
only need to bound $\max_{j=1,\ldots,N-1}A_j$ away from zero but that we require a suitable quantitative estimate of $\max_{j=1,\ldots,N-1}A_j$ from below.

Let us begin with the following definitions.
Let us call $\tilde{\Pi}=\mathrm{span}\{v_1,\ldots,v_{N-1}\}$.
For any $\Sigma\in\mathcal{A}_{obst}$, we denote with 
$P(\Sigma)$ its center of mass and $\tilde{\Pi}(\Sigma)=\tilde{\Pi}+P(\Sigma)$.
We define $\mathcal{A}_{sym}$ the set of $\Sigma\in \mathcal{A}_{obst}$ such that
$\Sigma$ is symmetric with respect to $\tilde{\Pi}(\Sigma)$.

We then define
the metric space
$$X=\{\Pi:\ \Pi\text{ is a hyperplane in }\mathbb{R}^N\text{ passing through }\overline{B_{R_2+1}}\}$$
with the distance
$$d(\Pi_1,\Pi_2)=
%d_H(\Sigma_1,\Sigma_2)+
d_H(\Pi_1\cap\overline{B_{2R_2+1}},\Pi_2\cap\overline{B_{2R_2+1}})\quad\text{for any }\Pi_1,\Pi_2\in X.
$$
Finally, we call $\mathcal{X}=\mathcal{A}_{obst}\times X$, with the standard metric of the product of two metric spaces, and
$\mathcal{Y}=\{(\Sigma,\tilde{\Pi}(\Sigma)):\ \Sigma\in\mathcal{A}_{sym}\}\subset\mathcal{X}.$

We have the following preliminary properties.

\begin{prop}
We have that $\Sigma\to P(\Sigma)$ is a Lipschitz continuous function on $\mathcal{A}_{obst}$ endowed with the Hausdorff distance, with a Lipschitz constant depending on the a priori data only.
Consequently $\mathcal{A}_{sym}$ is a closed subset of $\mathcal{A}_{obst}$ and
$\mathcal{Y}$ is closed in $\mathcal{X}$
\end{prop}

\proof{.} This is a straightforward consequence of 
\eqref{boundary}, \eqref{boundary33} and \eqref{distancesoppositebis}.\cvd

\begin{lem}\label{symmetriccase}
Let $\Sigma\in\mathcal{A}_{sym}$.
For simplicity, let us assume that $\tilde{\Pi}(\Sigma)=\{y_N=0\}$.
Then we call $G^{\pm}=\left\{y\in G:\ y_N\gtrless 0\right\}$ and we have that 
$G^{\pm}$ are Lipschitz domains with constants depending on the a priori data only.
\end{lem}

\proof{.}
The difficult part of the proof is to consider the points $z$ of $\partial G^{\pm}$
such that $z\in\partial\Sigma\cap \tilde{\Pi}(\Sigma)$

Let us consider a point $z\in\partial\Sigma\cap \tilde{\Pi}(\Sigma)$. By the Lipschitz properties of $\Sigma$, we have that there exists a given cone $\mathcal{C}$, with vertex in $0$, such that, for any
$y\in\partial\Sigma\cap B_{r/2}(z)$, $y+\mathcal{C}\subset G$. 
Since $\Sigma$ is symmetric with respect to $\tilde{\Pi}(\Sigma)$, we also have that
$y+\tilde{T}(\mathcal{C})\subset G$, $\tilde{T}$ being the reflection in $\tilde{\Pi}(\Sigma)$. Hence it is not difficult to show that there exists a cone $\mathcal{C}_1$, with vertex in $0$ and symmetric with respect to $\tilde{\Pi}(\Sigma)$, such that,
for any
$y\in\partial\Sigma\cap B_{\tilde{r}}(z)$, $y+\mathcal{C}_1\subset G$. We notice that
$\tilde{r}>0$ and
 the amplitude of the cone $\mathcal{C}_1$ depends on $r$ and $L$ only.

Therefore, for any point $z\in\partial\Sigma\cap \tilde{\Pi}(\Sigma)$, locally in $B_{\tilde{r}_1}(z)$,
$\partial\Sigma$ is the graph of a Lipschitz function, with Lipschitz constant bounded by $\tilde{L}_1$, with respect to a Cartesian coordinate system such that $e_N\in \tilde{\Pi}$. Hence it is not difficult to show that, locally in $B_{\tilde{r}_2}(z)$ and with respect to a different Cartesian coordinate system,
$\partial G^+$, and by symmetry $G^-$ as well, is the graph of a Lipschitz function, with Lipschitz constant bounded by $\tilde{L}_2$.
Here $\tilde{r}_i$ and $\tilde{L}_i$, $i=1,2$, are positive constants depending on $r$ and $L$ only.  

The proof now can be easily concluded.\cvd

\bigskip

In order to obtain the required lower bound on $\max_{j=1,\ldots,N-1}A_j$, we distinguish between two cases. The good one is when either $\Sigma$ is not close to $\mathcal{A}_{sym}$ or, if it is, the hyperplane $\Pi_M$ is not close to $\tilde{\Pi}(\Sigma)$. The bad one is when $\Sigma$ is close to $\mathcal{A}_{sym}$ and the hyperplane $\Pi_M$ is close to $\tilde{\Pi}(\Sigma)$.

In the next proposition we deal with the good case, in the sequel of the proof we shall show that, by a suitable modification of our geometric construction, the bad case actually never occurs.

Let us consider 
the map
$$\mathcal{X}\ni(\Sigma,\Pi)\mapsto f(\Sigma,\Pi)=\max_{j=1,\ldots,N-1}
\left(\max_{\tilde{z}\in (\overline{B_{2R_2+3}}\backslash B_{2R_2+2})
\cap \Pi}|\nabla u_j(\tilde{z})\cdot \nu|\right),$$
where $\nu=\nu_{\Pi}$ is the normal to $\Pi$ and, for any $j=1,\ldots,N-1$, $u_j$ is the solution to the direct scattering problem \eqref{uscateq} with boundary condition \eqref{soundhard}
and incident field $u^i(x)=\rme^{\rmi k x\cdot v_j}$, $x\in\mathbb{R}^N$. Then the following result holds.

\begin{prop}\label{goodcaseprop}
Let us fix a positive constant $\tilde{c}$. For any $a>0$, let us consider the following subset $\mathcal{X}_a$ of $\mathcal{X}$.
We say that $(\Sigma,\Pi)\in\mathcal{X}$ belongs to $\mathcal{X}_a$ if there exists
$\hat{\Sigma}\in\mathcal{A}_{sym}$ such that $d_H(\Sigma,\hat{\Sigma})<a$ and $d(\Pi,\tilde{\Pi}(\Sigma))<\tilde{c}a$.
 
Then there exists a positive constant $\hat{a}_0$, depending on the a priori data, on
$\tilde{c}$, on
 $a$ and on $\{v_1,\ldots,v_{N-1}\}$ only, such that
$$\min
\left\{f(\Sigma,\Pi):\ (\Sigma,\Pi)\in \mathcal{X}\backslash \mathcal{X}_a\right\}
\geq \hat{a}_0.$$
\end{prop}

\begin{oss}
In the previous proposition, if $N=2$ the result does not depend on the direction $v$. If $N=3$, the dependence on $v_1$ and $v_2$ is only through the constant
$b_0=|v_1\cdot v_2|<1$. We also notice that $\hat{a}_0$ does not depend on $h$.
\end{oss}

\proof{.} 
We observe that, by the stability result for the direct scattering problem with respect to sound-hard scatterers $\Sigma$ proved in \cite{Men-Ron}, such a map $f$ is continuous on $\mathcal{X}$.

If $f(\Sigma,\Pi)=0$, then the unit normal to $\Pi$ is orthogonal to $v_j$, for any $j=1,\ldots,N-1$, and $\Sigma$ is symmetric with respect to $\Pi$, that is $\Pi=\tilde{\Pi}(\Sigma)$ and $(\Sigma,\Pi)\in\mathcal{Y}$.

Then the proof immediately follows by the fact that $\mathcal{X}\backslash \mathcal{X}_a$ is closed and obviously does not contain any point of $\mathcal{Y}$.\cvd

\bigskip

Up to now, we are able to prove a stability result if either the assumptions of
Proposition~\ref{easycaseprop} are satisfied or, otherwise, if $(\Sigma,\Pi_M)
\in \mathcal{X}\backslash \mathcal{X}_a$ for a suitable $a>0$.
In both cases we use the same computation as in the $N$ measurements case, with $a_0$ replaced by $\tilde{a}_0$ and $\hat{a}_0$, respectively. We notice that in this second case $\hat{a}_0$ depends on $a$.

Therefore our strategy is now the following. We choose a suitable value of $a$ and we construct a modified regular chain for $\Sigma$ as in the general geometric construction such that for any possible reflection hyperplane $\Pi_n$, $n=1,\ldots,M$, (including the first one!) we have that $(\Sigma,\Pi_n)\in \mathcal{X}\backslash \mathcal{X}_a$.
As we shall see, actually the first reflection hyperplane is the one that presents the greatest difficulties.

We notice that, so far, all our arguments work for any dimension $N\geq 2$. However the construction of such a modified regular chain presents some technical challenges,
in particular for the proof of Lemma~\ref{startingpointlemma} below.
 Therefore in the sequel we limit ourselves to the space dimension $N=3$ and we notice that when the space dimension is $N=2$ the result may be proved along the same lines.

A crucial remark is that, unfortunately, we are not able to choose $a$ independently of $h$. This is the main reason why we lose the precise dependence of our stability result on the size parameter $h$, that we instead have in the sound-soft case or in the sound-hard case with $N$ measurements.

The following two technical lemmas shall be needed.

\begin{lem}\label{G+Lipprop}
Let $N= 3$ and $h>0$.

There exist positive constants $\tilde{c}_0$, $\tilde{c}_1$, $\tilde{c}_2\leq 1$  and $\tilde{L}\geq L$, depending on the a priori data only, such that the following holds.

Let $a=\tilde{c}_0h$
and let $\Sigma\in \mathcal{A}^h_{obst}$ satisfy the following. We assume that there exists
$\hat{\Sigma}\in\mathcal{A}_{sym}$ such that $d_H(\Sigma,\hat{\Sigma})\leq a$.

For simplicity, let us assume that $\tilde{\Pi}(\Sigma)=\{y_N=0\}$.
Then, for any $\tilde{c}$, $0\leq \tilde{c}\leq \tilde{c}_1$, if
we call $G^{\pm}_a=\{y\in G:\ y_N\gtrless \pm\tilde{c}a\}$, we have that 
$G^{\pm}_a$ are Lipschitz domains with constants $\tilde{r}=\tilde{c}_2h$ and $\tilde{L}$.
\end{lem}

\proof{.} This is an extension of Lemma~\ref{symmetriccase}, which can be proved by exploiting
\cite[Proposition~6.1]{Ron08}.\cvd

\begin{lem}\label{startingpointlemma}
Let $N= 3$, $h>0$ and $\tilde{c}_0$ and $\tilde{c}_1$ as in Lemma~\textnormal{\ref{G+Lipprop}}.
Let $\Sigma$, $\Sigma'$ belong to $\mathcal{A}^h_{obst}$ and let $d$ be defined as in \eqref{ddefin}.
Let $x_1\in \partial\Sigma'\backslash\Sigma$ be such that
$d=\mathrm{dist}(x_1,\partial\Sigma)=\mathrm{dist}(x_1,\Sigma)$. 

There exist positive constants $\tilde{c}_3\leq 1$, $\tilde{c}_4\leq \tilde{c}_1$ and $K_1\leq 1$, depending on the a priori data only, such that the following holds.

Let $a=\tilde{c}_0h$ and $\tilde{c}=\tilde{c}_4$.
Let us assume that there exists
$\hat{\Sigma}\in\mathcal{A}_{sym}$ such that $d_H(\Sigma,\hat{\Sigma})\leq a$.
Let us call $G^{\pm}_a=\{y\in G:\ y_N\gtrless \pm\tilde{c}a\}$,
assuming that $\tilde{\Pi}(\Sigma)=\{y_N=0\}$.

If $d\leq \tilde{c}_3h$, then, up to swapping the role of $G^+_a$ and $G^-_a$, there exists $\tilde{x}_1\in\partial\Sigma'\backslash\Sigma$
such that $\tilde{x}_1\in G^+_a$ and 
\begin{equation}\label{newstartingpoint}
\mathrm{dist}(\tilde{x}_1,\partial G^+_a)\geq K_1d^3.
\end{equation}
\end{lem}

\proof{.} This is a straightforward consequence of \cite[Proposition~6.2]{Ron08}.\cvd

\bigskip

We are now in the position of concluding the proof of the $N-1$ measurements case

\proof{ of Theorem~\ref{mainteoN-1}.} Without loss of generality we can assume that
$h\leq \min\{r,1\}$.

Let us assume, for the time being, that $d\leq \tilde{c}_3h\leq h$.
Let us set $a=\tilde{c}_0h$, $\tilde{c}_0$ as in Lemma~\ref{G+Lipprop},
and $\tilde{c}=\tilde{c}_4$
as in Lemma~\ref{startingpointlemma}.

Then we distinguish between two cases. If there does not exist any
$\hat{\Sigma}\in\mathcal{A}_{sym}$ such that $d_H(\Sigma,\hat{\Sigma})\leq a$,
then we conclude using the geometric construction of Section~\ref{geomconstrsec} and the arguments used for the proof of the $N$ measurements case. Here we use either Proposition~\ref{easycaseprop}, replacing $a_0$ with $\tilde{a}_0$, or Proposition~\ref{goodcaseprop}, with $\tilde{c}=\tilde{c}_4$ as in Lemma~\ref{startingpointlemma} and replacing $a_0$ with $\hat{a}_0$. We have to notice that $\hat{a}_0$ here depends on $a$ thus on $h$.

Otherwise, let us assume that there does exist
$\hat{\Sigma}\in\mathcal{A}_{sym}$ such that $d_H(\Sigma,\hat{\Sigma})\leq a$.
Then we use the geometric construction and estimates of
Section~\ref{geomconstrsec} with the following differences. We replace $x_1$ with $\tilde{x}_1$ and $G$ with $G^+_a$, $\tilde{x}_1$ and $G^+_a$ as in Lemma~\ref{startingpointlemma}. Using Lemma~\ref{G+Lipprop}, we further replace $r$ and $L$
with $\tilde{c}_2h$ and $\tilde{L}$, respectively. Finally, using
\eqref{newstartingpoint}, we replace $d$ with $K_1d^3\leq d$.
Then we can repeat the previous argument using either 
Proposition~\ref{easycaseprop} or Proposition~\ref{goodcaseprop}.
In fact any possible reflection point belongs to $G^+_a$,
therefore any reflection plane is far enough from $\tilde{\Pi}(\Sigma)$.

We conclude that, for any $\varepsilon$, $0<\varepsilon<1/(2\rme)$, 
provided $d\leq \tilde{c}_3h$,
we have
$$
K_1d^3 \leq 2\rme R(\eta(\varepsilon))^{3C}
$$
for some constant $C$ depending on the a priori data, on $b_0$ and on $h$ as well.
Therefore
\begin{equation}\label{almostdone}
d\leq A_1(\eta(\varepsilon))^C
\end{equation}
where $A_1$ depends on the a priori data only and $C$ depends on the a priori data, on $b_0$ and on $h$.

Finally, we need to drop the assumption that $d\leq \tilde{c}_3h$.
We claim that there exists $\hat{\varepsilon}_1(h)$, $0<\hat{\varepsilon}_1(h)\leq 1/(2\rme)$, depending on the a priori data and on $h$ only, such that
\begin{multline}\label{finalclaim}
\hat{\varepsilon}_1(h)<\\
\inf\left\{
\|u-u'\|_{L^{\infty}(B_{\tilde{\rho}}(x_0))}:\ v\in\mathbb{S}^{N-1},\ \Sigma,\Sigma'\in \mathcal{A}^h_{obst}\text{ such that }d_H(\Sigma,\Sigma')
\geq C_1\tilde{c}_3h   
\right\}
\end{multline}
where $C_1$ is as in 
\eqref{distancesoppositebis}. If this is true, then obviously we obtain that if $\varepsilon\leq \hat{\varepsilon}_1(h)$ then 
$d_H(\Sigma,\Sigma')< C_1\tilde{c}_3h$ that is
$d\leq d_H(\Sigma,\Sigma')/C_1<\tilde{c}_3h$ and the proof would be concluded.

Therefore we just need to prove the claim in \eqref{finalclaim}. It is not difficult to show that the infimum on the right hand side is actually a minimum. Again it is enough to use the stability result of the direct scattering problem with respect to the variation of sound-hard scatterers proved in \cite{Men-Ron}. Finally, if such a minimum were zero we would contradict the uniqueness result for the determination of a sound-hard obstacle by a single scattering measurement proved in \cite{Els-Yam1,Els-Yam2}.\cvd

\subsection{The single measurement case}\label{1meassubs}

We restrict here to $N=3$, since $N=2$ is clearly covered by the previous section. We consider the assumptions and notation of Theorem~\ref{mainteo1} to hold.

The main technical difficulty we have to tackle if we have only one measurement, compared to the two measurements case, is that we may have several planes whose normal is orthogonal to $v$ with respect to which $\Sigma$ might be symmetric. As we discussed in the previous subsection, using two measurements with two directions of propagation $v_1$ and $v_2$ allows us to consider only one possible symmetry plane for $\Sigma$.

We begin with the following definition. Here we call $\mathcal{A}_{sym}$, respectively $\mathcal{A}^h_{sym}$, the set of $\Sigma$ belonging to $\mathcal{A}_{obst}$,
respectively $\mathcal{A}^h_{obst}$, such that $\Sigma$ is symmetric with respect to at least one plane whose normal is orthogonal to the incident direction of propagation $v$.
Moreover, for any $\Sigma\in \mathcal{A}^h_{obst}$ we call $n(\Sigma)$ the number of planes whose normal is orthogonal to $v$ with respect to which $\Sigma$ is symmetric. Notice that $n(\Sigma)$ is always a nonnegative integer that, obviously, could also be zero.
In other words, $\mathcal{A}^h_{sym}$ is the set of $\Sigma\in \mathcal{A}^h_{obst}$
such that $n(\Sigma)>0$.

We shall use the following notation. For any $\Sigma\in\mathcal{A}^h_{sym}$ we call $\Pi_i(\Sigma)$, $i=1,\ldots,n(\Sigma)$, the planes whose normal is orthogonal to $v$ with respect to which $\Sigma$ is symmetric. Correspondingly, we define $\nu_i(\Sigma)$, $i=1,\ldots,n(\Sigma)$, their corresponding unit normals, noticing that they all belong to the plane that is orthogonal to $v$.

We have the following properties whose proof is elementary and will be omitted.

\begin{prop}\label{symmprop}
There exists an integer $M=M(h)$, depending on the a priori data and on $h$ only, such that $n(\Sigma)\leq M$ for any $\Sigma\in\mathcal{A}^h_{sym}$. As a consequence, there exists a constant $\alpha=\alpha(h)$, $0<\alpha<\pi/2$, such that, for any $\Sigma\in\mathcal{A}^h_{sym}$, the angle between $\nu_i(\Sigma)$ and $\nu_j(\Sigma)$, with $i\neq j$, is bounded from below by $\alpha$.

We consider $\mathcal{A}_{obst}$ endowed with the Hausdorff distance.
Then the map
$\mathcal{A}^h_{obst}\ni\Sigma\mapsto n(\Sigma)$ is upper semicontinuous. Consequently, $\mathcal{A}^h_{sym}$ is a compact subset of $\mathcal{A}^h_{obst}$.
\end{prop}

Let us define, for any $n=1,\ldots,M(h)$, $\mathcal{A}^h_{sym,n}$ as the set of $\Sigma\in \mathcal{A}^h_{sym}$ such that $n(\Sigma)=n$.

The crucial difference with respect to the $2$ measurements case is that we need to define 
the set $\mathcal{X}_a$, for a positive constant $a$, in a rather more involved way. Next we describe such a construction, for any positive $a$ and for a fixed constant $\tilde{c}$ to be decided later.

Given $M=M(h)$ we begin in the following way. For any $\hat{\Sigma}\in \mathcal{A}^h_{sym,M}$, we find $r(\hat{\Sigma})$, $0<r(\hat{\Sigma})\leq a$, such that for any $\Sigma\in \mathcal{A}^h_{sym}$ with $d_{H}(\Sigma,\hat{\Sigma})\leq r(\hat{\Sigma})$ the following holds. For any $i=1,\ldots,n(\Sigma)$ there exists $j\in \{1,\ldots,n(\hat{\Sigma})\}$ such that $d(\Pi_i(\Sigma),\Pi_j(\hat{\Sigma}))\leq \tilde{c}a/4$.

Then, by compactness, we have that
$$\mathcal{A}^h_{sym,M}\subset \bigcup_{j=1}^{m_M}B_{r(\hat{\Sigma}_j)/2}(\hat{\Sigma}_j)=A_M^h$$
where, for any $j=1,\ldots,m_M$, $\hat{\Sigma}_j\in \mathcal{A}^h_{sym,M}$ and
$r(\hat{\Sigma}_j)\leq r(\hat{\Sigma}_{j-1})$. Here by convention we set $r(\hat{\Sigma}_0)=a$.

Then we consider $\mathcal{A}^h_{sym,M-1}\backslash A_M^h$, which is again a compact set. We consider a similar construction as before. Namely,
for any $\hat{\Sigma}\in \mathcal{A}^h_{sym,M-1}\backslash A_M^h$, we find $r(\hat{\Sigma})$, $0<r(\hat{\Sigma})\leq r(\hat{\Sigma}_{m_M})$, such that for any $\Sigma\in \mathcal{A}^h_{sym}$ with $d_{H}(\Sigma,\hat{\Sigma})\leq r(\hat{\Sigma})$ the following holds. For any $i=1,\ldots,n(\Sigma)$ there exists $j\in \{1,\ldots,n(\hat{\Sigma})\}$ such that $d(\Pi_i(\Sigma),\Pi_j(\hat{\Sigma}))\leq \tilde{c}a/4$. Then, by compactness, we have that
$$\mathcal{A}^h_{sym,M-1}\backslash A_M^h\subset \bigcup_{j=m_M+1}^{m_{M-1}}B_{r(\hat{\Sigma}_j)/2}(\hat{\Sigma}_j)=A_{M-1}^h$$
where, for any $j=m_M+1,\ldots,m_{M-1}$, $\hat{\Sigma}_j\in \mathcal{A}^h_{sym,M-1}\backslash A_M^h$ and
$r(\hat{\Sigma}_j)\leq r(\hat{\Sigma}_{j-1})$.

We proceed in a completely analogous way until we have that
$$\mathcal{A}^h_{sym}\subset \bigcup_{j=1}^{m_1}B_{r(\hat{\Sigma}_j)/2}(\hat{\Sigma}_j)=\bigcup_{l=1}^M A_l^h=A^h.$$
For any $l=2,\ldots,M$, and any
$j=m_l+1,\ldots,m_{l-1}$, $\hat{\Sigma}_j\in \mathcal{A}^h_{sym,l-1}\backslash \left(\bigcup_{i=l}^M A_i^h\right)$ and
$r(\hat{\Sigma}_j)\leq r(\hat{\Sigma}_{j-1})$.

We call $\mathcal{X}^h$ the subset of $(\Sigma,\Pi)\in\mathcal{X}$ such that
 $\Sigma\in \mathcal{A}^h_{obst}$.
We also call $\mathcal{Y}^h$ the subset of $\mathcal{X}^h$ defined as follows
$$\mathcal{Y}^h=\{(\Sigma,\Pi)\in\mathcal{X}^h:\ \Sigma\in \mathcal{A}^h_{sym}\text{ and }\Pi=\Pi_i(\Sigma)\text{ for some }i\in\{1,\ldots,n(\Sigma)\}\}.$$

Then we define $\mathcal{X}_a^h$ the subset of $\mathcal{X}^h$ with the following properties. We say that $(\Sigma,\Pi)\not\in \mathcal{X}_a^h$ either if $\Sigma\not\in A^h$ or if 
$\Sigma\in A^h$ and $d(\Pi,\hat{\Pi})\geq \tilde{c}a/2$ for any plane $\hat{\Pi}$ such that 
$\hat{\Pi}=\Pi_i(\hat{\Sigma}_j)$ for some $j\in\{1,\ldots,m_1\}$ such that
$d_H(\Sigma,\hat{\Sigma}_j)<r(\hat{\Sigma}_j)/2$ and for some $i\in\{1,\ldots,n(\hat{\Sigma}_j)\}$.

It is an easy remark that, for any $a>0$, $\mathcal{Y}^h\subset \mathcal{X}_a^h$
and $\mathcal{X}^h\backslash\mathcal{X}_a^h$ is closed. 
Let us consider 
the map
$$\mathcal{X}^h\ni(\Sigma,\Pi)\mapsto f(\Sigma,\Pi)=
\max_{\tilde{z}\in (\overline{B_{2R_2+3}}\backslash B_{2R_2+2})
\cap \Pi}|\nabla u(\tilde{z})\cdot \nu|,$$
where $\nu=\nu_{\Pi}$ is the normal to $\Pi$ and $u$ is the solution to the direct scattering problem \eqref{uscateq} with boundary condition \eqref{soundhard}
and incident field $u^i(x)=\rme^{\rmi k x\cdot v}$, $x\in\mathbb{R}^N$. Hence, arguing as in the proof of Proposition~\ref{goodcaseprop}, we can obtain the following result.

\begin{prop}\label{goodcaseprop1}
Let us fix a positive constant $\tilde{c}$. For any $a>0$, we define the subset $\mathcal{X}_a^h$ of $\mathcal{X}$ as before.
 
Then there exists a positive constant $\hat{a}_0$, depending on the a priori data, on
$\tilde{c}$,  on
 $a$ and on $h$ only, such that
$$\min
\left\{f(\Sigma,\Pi):\ (\Sigma,\Pi)\in \mathcal{X}^h\backslash \mathcal{X}^h_a\right\}
\geq \hat{a}_0.$$
\end{prop}

We now consider the corresponding results to Lemmas~\ref{G+Lipprop} and \ref{startingpointlemma}. We need the following notation, recalling that positive constants $a$ and $\tilde{c}$ are fixed.
For any $\Sigma\in A^h$ we choose $\hat{\Sigma}(\Sigma)$ as the first 
$\hat{\Sigma}_j$, $j\in\{1,\ldots,m_1\}$, such that $\Sigma\in B_{r(\hat{\Sigma}_j)/2}(\hat{\Sigma}_j)$. For any $i=1,\ldots,n(\hat{\Sigma}(\Sigma))$, we define the infinite strips
$$S_i=\Pi_i(\hat{\Sigma}(\Sigma))+\{ca\nu_i(\hat{\Sigma}(\Sigma)):\ |c|\leq\tilde{c}\}.$$
We notice that $\mathbb{R}^3\backslash (\bigcup_{i=1}^{n(\hat{\Sigma}(\Sigma))}S_i)$
consists of $2n(\hat{\Sigma}(\Sigma))$ different connected open sectors that we shall call $G_a^j$,
$j=1,\ldots,2n(\hat{\Sigma}(\Sigma))$.

Then, with the notation introduced above,  the following important results hold.

\begin{lem}\label{G+Lipprop1}
Let $N= 3$ and $h>0$. 

There exist positive constants $\tilde{c}_0$, depending on the a priori data only,
and $\tilde{c}_1$, $\tilde{c}_2\leq 1$ and $\tilde{L}\geq L$, depending on the a priori data and on $h$ only, such that the following holds.

Let $a=\tilde{c}_0h$
and let $\Sigma\in A^h$.
Then, for any $\tilde{c}$, $0\leq \tilde{c}\leq \tilde{c}_1$, we have that, for any $j=1,\ldots,2n(\hat{\Sigma}(\Sigma))$, 
$G_a^j\backslash \Sigma$ is a Lipschitz domain with constants $\tilde{r}=\tilde{c}_2h$ and $\tilde{L}$.
\end{lem}

\proof{.} If $\hat{\Sigma}(\Sigma)\in \mathcal{A}^h_{sym,1}$, then the result is contained in Lemma~\ref{G+Lipprop}. Therefore, without loss of generality we assume that
$\hat{\Sigma}(\Sigma)\in \mathcal{A}^h_{sym,n}$ for some $n\geq 2$.

We begin by proving the following claim, which is the corresponding result to
Lemma~\ref{symmetriccase}.
 We fix an arbitrary $\hat{\Sigma}\in\mathcal{A}^h_{sym,n}$ with $n\geq 2$.
Then we call $G^j$, $j=1,\ldots,2n$, the connected components of $\mathbb{R}^3\backslash(\bigcup_{i=1}^n\Pi_i(\hat{\Sigma}))$. We claim that, for any $j=1,\ldots,2n$, $G_j\backslash\hat{\Sigma}$ is a Lipschitz domain with constants $\tilde{r}_1$ and $\tilde{L}_1$ depending on the a priori data and on $h$ only.

We deal with the points $z$ belonging to $\partial\hat{\Sigma}$ and $\Pi_i(\hat{\Sigma})$ for some $i=1,\ldots,n$. Let $P(\hat{\Sigma})$ be the center of mass of $\hat{\Sigma}$ and let $l$ be the line defined as follows
$$l=\{x\in\mathbb{R}^3: x=P(\hat{\Sigma})+rv,\ r\in\mathbb{R}\}.$$
It is obvious that $l=\Pi_i(\hat{\Sigma})\cap\Pi_j(\hat{\Sigma})$ for any $i\neq j$.

If we have a point $z$ belonging to $\partial\hat{\Sigma}$ and $\Pi_i(\hat{\Sigma})$ for some $i=1,\ldots,n$, which is far enough from $l$, we can treat it exactly as in the proof of Lemma~\ref{symmetriccase}. Therefore the most delicate case is the one in which $z\in\partial\hat{\Sigma}\cap l$. However, following the kind of reasonings used in the proof of Lemma~\ref{symmetriccase}, it is not difficult to show that, locally in $B_{\tilde{r}_2}(z)$,
$\partial\Sigma$ is the graph of a Lipschitz function, with Lipschitz constant bounded by $\tilde{L}_2$, with respect to a Cartesian coordinate system such that $e_3$ is parallel to $v$, with $\tilde{r}_2$ and $\tilde{L}_2$ depending on the a priori data only.

Then the claim easily follows, with the dependence of $\tilde{r}_1$ and $\tilde{L}_1$ on $h$ essentially given by the angle $\alpha(h)$.

The proof of the proposition can be concluded by
using the claim, arguments similar to the ones used to prove the claim, and
\cite[Proposition~6.1]{Ron08}.\cvd

\bigskip

We notice that the difference with respect to Lemma~\ref{G+Lipprop} is that now $\tilde{c}_1$, $\tilde{c}_2$ and $\tilde{L}$ depend on $h$ as well.

\begin{lem}\label{startingpointlemma1}
Let $N= 3$, $h>0$ and $\tilde{c}_0$ and $\tilde{c}_1$ as in Lemma~\textnormal{\ref{G+Lipprop1}}.
Let $\Sigma$, $\Sigma'$ belong to $\mathcal{A}^h_{obst}$ and let $d$ be defined as in \eqref{ddefin}.
Let $x_1\in \partial\Sigma'\backslash\Sigma$ be such that
$d=\mathrm{dist}(x_1,\partial\Sigma)=\mathrm{dist}(x_1,\Sigma)$. 

There exist positive constants $\tilde{c}_3\leq 1$, $\tilde{c}_4\leq \tilde{c}_1$,
depending on the a priori data and on $h$ only,
and $K_1\leq 1$, depending on the a priori data only, such that the following holds.

Let $a=\tilde{c}_0h$ and $\tilde{c}=\tilde{c}_4$.
Let us assume that $\Sigma\in A^h$.

If $d\leq \tilde{c}_3h$, then there exist $\tilde{x}_1\in\partial\Sigma'\backslash\Sigma$
and
 $j\in\{1,\ldots,2n(\hat{\Sigma}(\Sigma))\}$ such that
$\tilde{x}_1\in G^j_a\backslash \Sigma$ and 
\begin{equation}\label{newstartingpoint1}
\mathrm{dist}(\tilde{x}_1,\partial (G^j_a\backslash \Sigma))\geq K_1d^3.
\end{equation}
\end{lem}

\proof{.} The result is a rather straightforward consequence of \cite[Proposition~6.2]{Ron08}.\cvd

\bigskip

Again, it is important to remark that 
the difference with respect to Lemma~\ref{startingpointlemma} is that now $\tilde{c}_3$ and $\tilde{c}_4$ depend on $h$ too.

We are now in the position of proving our stability result with one measurement.

\proof{ of Theorem~\ref{mainteo1}.} The proof follows the same arguments of the proof of Theorem~\ref{mainteoN-1}, replacing Proposition~\ref{goodcaseprop}, Lemmas~\ref{G+Lipprop} and \ref{startingpointlemma} with Proposition~\ref{goodcaseprop1}, Lemmas~\ref{G+Lipprop1} and \ref{startingpointlemma1}, respectively. We point out the modification that we need to adopt in this case.

Without loss of generality we can assume that
$h\leq \min\{r,1\}$.
Let us assume, for the time being, that $d\leq \tilde{c}_3h\leq h$.
Let us set $a=\tilde{c}_0h$, $\tilde{c}_0$ as in Lemma~\ref{G+Lipprop1},
and $\tilde{c}=\tilde{c}_4$
as in Lemma~\ref{startingpointlemma1}.

We distinguish two cases. If $\Sigma\not\in A^h$, then we conclude using Proposition~\ref{easycaseprop} and Proposition~\ref{goodcaseprop1} with $\tilde{c}=\tilde{c}_4$ as in Lemma~\ref{startingpointlemma1}.

If instead $\Sigma\in A^h$, we replace $x_1$ with $\tilde{x}_1$ and $G$ with $G^j_a\backslash \Sigma$, $\tilde{x}_1$ and
$G^j_a\backslash \Sigma$ as in Lemma~\ref{startingpointlemma1}.
Notice that in this case, any possible reflection point belongs to $G^j_a\backslash \Sigma$, therefore any reflection plane $\Pi$ is far from any $\Pi_i(\hat{\Sigma}(\Sigma))$, $i=1,\ldots,n(\hat{\Sigma}(\Sigma))$, at least $\tilde{c}a$. On the other hand, recalling the construction of $A^h$ and how we choose $\hat{\Sigma}(\Sigma)$, for any $j\in\{1,\ldots,m_1\}$, if
$d_H(\Sigma,\hat{\Sigma}_j)<r(\hat{\Sigma}_j)/2$, we have that
$d_H(\hat{\Sigma}_j,\hat{\Sigma}(\Sigma))<r(\hat{\Sigma}(\Sigma))$. We conclude that
$d(\Pi,\hat{\Pi})\geq 3\tilde{c}a/4$ for any plane $\hat{\Pi}$ such that 
$\hat{\Pi}=\Pi_i(\hat{\Sigma}_j)$ for some $j\in\{1,\ldots,m_1\}$ such that
$d_H(\Sigma,\hat{\Sigma}_j)<r(\hat{\Sigma}_j)/2$ and for some $i\in\{1,\ldots,n(\hat{\Sigma}_j)\}$. That is $(\Sigma,\Pi)\in\mathcal{X}^h\backslash\mathcal{X}_a^h$.

The rest of the argument is the same. However, we notice that, since the domains $G^j_a\backslash \Sigma$ used in Lemmas~\ref{G+Lipprop1} and \ref{startingpointlemma1} are Lipschitz with constants both depending on $h$, the dependence of the stability result on $h$ may be worse than the one in the $2$ measurements case.\cvd


\begin{thebibliography}{99}

%\bibitem{Ada}
%R.~A.~Adams,
%\emph{Sobolev Spaces},
%Academic Press, New York, 1975.

\bibitem{Ale-Ron}
G.~Alessandrini and L.~Rondi,
\emph{Determining a sound-soft polyhedral scatterer by a single far-field
measurement},
Proc. Amer. Math. Soc. \textbf{133} (2005) 1685--1691.

\bibitem{Bru}
R.~Brummelhuis,
\emph{Three-spheres theorem for second order elliptic equations},
J. Analyse Math. \textbf{65} (1995) 179--206.

\bibitem{Che-Yam}
J.~Cheng and M.~Yamamoto,
\emph{Uniqueness in an inverse scattering problem within non-trapping polygonal obstacles with at most two incoming waves},
Inverse Problems \textbf{19} (2003) 1361--1384.

\bibitem{Col-Kre98}
D.~Colton and R.~Kress,
\emph{Inverse Acoustic and Electromagnetic Scattering Theory},
Springer-Verlag, Berlin Heidelberg New York, 1998.

\bibitem{Col-Sle}
D.~Colton and B.~D.~Sleeman,
\emph{Uniqueness theorems for the inverse problem of acoustic scattering},
IMA J. Appl. Math. \textbf{31} (1983) 253--259. 

\bibitem{Els-Yam1}
J.~Elschner and M.~Yamamoto,
\emph{Uniqueness in determining polygonal sound-hard obstacles with a single incoming wave},
Inverse Problems \textbf{22} (2006) 355--364.

\bibitem{Els-Yam2}
J.~Elschner and M.~Yamamoto,
\emph{Uniqueness in determining polyhedral sound-hard obstacles with a single incoming wave},
Inverse Problems \textbf{24} (2008) 035004 (7pp).

\bibitem{Gia}
A.~Giacomini,
\emph{A stability result for Neumann problems in dimension} $N\geq 3$,
J. Convex Anal. \textbf{11} (2004) 41--58.

\bibitem{Isak90}
V.~Isakov,
\emph{On uniqueness in the inverse transmission scattering problem},
Comm. Partial Differential Equations \textbf{15} (1990) 1565--1587.

\bibitem{Isak92}
V.~Isakov,
\emph{Stability estimates for obstacles in inverse scattering},
J. Comp. Appl. Math. \textbf{42} (1992) 79--88.

\bibitem{Isak93}
V.~Isakov,
\emph{New stability results for soft obstacles in inverse scattering},
Inverse Problems \textbf{9} (1993) 535--543.

\bibitem{Isak06}
V.~Isakov,
\emph{Inverse Problems for Partial Differential Equations}, 2nd edition,
Springer-Verlag, New York, 2006.

\bibitem{Kir-Kre}
A.~Kirsch and R.~Kress,
\emph{Uniqueness in inverse obstacle scattering},
Inverse Problems \textbf{9} (1993) 285--299.

\bibitem{L-P}
P.~D.~Lax and R.~S.~Phillips,
\emph{Scattering Theory}, Academic Press,  New York London, 1967.

\bibitem{Liu-Zou}
H.~Liu and J.~Zou,
\emph{Uniqueness in an inverse acoustic obstacle scattering problem for both sound-hard and sound-soft polyhedral scatterers},
Inverse Problems \textbf{22} (2006)
515--524.

\bibitem{Liu-Zou3}
H.~Liu and J.~Zou,
\emph{On unique determination of partially coated polyhedral scatterers
with far field measurements},
Inverse Problems \textbf{23} (2007) 297--308.

\bibitem{Men-Ron}
G.~Menegatti and L.~Rondi,
\emph{Stability for the acoustic scattering problem for sound-hard scatterers},
Inverse Probl. Imaging \textbf{7} (2013) 1307--1329.

%\bibitem{Pott}
%R.~Potthast,
%\emph{On a concept of uniqueness in inverse scattering for a finite number of incident waves},
%SIAM J. Appl. Math. \textbf{58} (1998) 666--682. 

\bibitem{Ron03}
L.~Rondi,
\emph{Unique determination of non-smooth sound-soft scatterers by finite\-ly many
far-field measurements},
Indiana Univ. Math. J. \textbf{52} (2003) 1631--1662.

\bibitem{Ron07}
L.~Rondi, \emph{A variational approach to the reconstruction of cracks by boundary measurements}, J. Math. Pures Appl. (9) \textbf{87} (2007) 324--342.

\bibitem{Ron08}
L.~Rondi,
\emph{Stable determination of sound-soft polyhedral scatterers by a single measurement},
Indiana Univ. Math. J. \textbf{57} (2008) 1377--1408.

\bibitem{Ron-Sin}
L.~Rondi and M.~Sini,
\emph{Stable determination of a scattered wave from its far-field
pattern}: \emph{the high frequency asymptotics},
Arch. Ration. Mech. Anal.
\textbf{218} (2015) 1--54.

%\bibitem{Try}
%G.~N.~Trytten,
%\emph{Pointwise bounds for solutions of the Cauchy problem for elliptic
%equations},
%Arch. Rational Mech. Anal. \textbf{13} (1963) 222--244. 

\bibitem{Wil}
C.~H.~Wilcox,
\emph{Scattering Theory for the d'Alembert Equation in Exterior Domains},
Springer-Verlag, Berlin New York, 1975.

\end{thebibliography}
\end{document}